\documentclass[journal]{IEEEtran}
\usepackage{lipsum}
\usepackage{graphicx}
\usepackage{epstopdf}
\usepackage{braket,amsfonts}
\usepackage{amsopn}
\usepackage{caption}
\usepackage{amsthm}
\usepackage{amsmath}
\usepackage{makecell}
\usepackage{xcolor}
\usepackage{siunitx}
\usepackage{booktabs}
\usepackage{blindtext}
\usepackage{hyperref}
\usepackage{tabularx,booktabs}
\usepackage{caption}
\usepackage{subcaption}
\usepackage{tablefootnote}
\newcolumntype{C}{>{\centering\arraybackslash}X} 
\newcommand{\Z}{\mathbb{Z}}
\newcommand{\R}{\mathbb{R}}

\newtheorem{definition}{Definition}
\newtheorem{assumption}{Assumption}

\newtheorem{theorem}{Theorem}

\usepackage{caption}
\usepackage{cite}

\usepackage{amsmath,amssymb}

\usepackage{algorithmic, algorithm}

\usepackage{array}

\usepackage{url}

\allowdisplaybreaks

\begin{document}
\title{A Two-level ADMM Algorithm for AC OPF with {Global} Convergence Guarantees}

\author{Kaizhao~Sun
        and Xu Andy~Sun,~\IEEEmembership{Senior Menmber,~IEEE}
\thanks{The authors are with the H. Milton Stewart 
of Industrial and Systems Engineering, Georgia Institute of Technology, Atlanta,
GA, 30332 USA (email: \texttt{ksun46@gatech.edu}; \texttt{andy.sun@isye.gatech.edu}).}
\thanks{Manuscript received August 27, 2020; revised January 28, 2021; {accepted March 28, 2021.}}}

\maketitle

\begin{abstract}
This paper proposes a two-level distributed algorithmic framework for solving the AC optimal power flow (OPF) problem with convergence guarantees. The presence of highly nonconvex constraints in OPF poses significant challenges to distributed algorithms based on the alternating direction method of multipliers (ADMM). In particular, convergence is not provably guaranteed for nonconvex network optimization problems like AC OPF. In order to overcome this difficulty, we propose a new distributed reformulation for AC OPF and a two-level ADMM algorithm that goes beyond the standard framework of ADMM. We establish the global convergence and iteration complexity of the proposed algorithm under mild assumptions. Extensive numerical experiments over some largest test cases from NESTA and PGLib-OPF (up to 30,000-bus systems) demonstrate advantages of the proposed algorithm over existing ADMM variants in terms of convergence, scalability, and robustness. Moreover, under appropriate parallel implementation, the proposed algorithm exhibits fast convergence comparable to or even better than the state-of-the-art centralized solver.
\end{abstract}

\begin{IEEEkeywords}
Distributed optimization, optimal power flow, augmented Lagrangian method, alternating direction method of multipliers.
\end{IEEEkeywords}

%
\IEEEpeerreviewmaketitle

\section{Introduction}
The AC optimal power flow (OPF) is a basic building block in electric power grid operation and planning. It is a highly nonconvex optimization problem, due to nonlinear power flow equations, and is shown to be an NP-hard decision problem \cite{bienstock2015strong,lehmann2016ac}. Any computational method to be deployed in power system practice should meet the stringent requirement that the algorithm has guaranteed robust performance and requires minimal tuning and intervention in face of variations in system conditions. Moreover, to effectively coordinate multiple regions in a large power grid, distributed algorithms that do not require sharing critical private information between regions should be highly desirable. 
However, such a goal has remained challenging for solving large-scale AC-OPF problems. Many existing algorithms are only suited for centralized operation. Most distributed or decentralized algorithms for solving AC-OPF do not have guaranteed convergence performance and require extensive parameter tuning and experimentation. 

In this paper, we develop a distributed algorithm for solving large-scale AC OPF problems {to stationary points}, and the proposed algorithm is proven to have global convergence. 



\subsection{Literature Review}
The research community has extensively studied local nonlinear optimization methods such as the interior point methods e.g. \cite{jabr2002primal, zimmerman2011matpower}. Another line of research looks into convex relaxations of AC OPF and has drawn significant attentions in recent years. 
In particular, the semidefinite programming (SDP) relaxation is firstly applied to the OPF problem in \cite{bai2008semidefinite}, and sufficient conditions are studied to guarantee the exactness of SDP relaxations \cite{lavaei2014geometry}.
However, SDP suffers from expensive computation cost for large-scale problems, while the second-order cone programming (SOCP) relaxation initially proposed in \cite{jabr2006radial} offers a favorable alternative. The strong SOCP relaxation proposed in \cite{kocuk2016strong} is shown to be close to or dominates the SDP relaxation, while the computation time of SOCP can be orders of magnitude faster than SDP. However, obtaining a primal feasible solution is a common challenge facing these convex relaxation methods. 

The alternating direction method of multipliers (ADMM) offers a powerful framework for distributed computation. Sun et al. \cite{sun2013fully} applied ADMM to decompose the computation down to each individual bus, and observe the convergence is sensitive to initial conditions. Erseghe \cite{erseghe2014distributed, erseghe2015distributed} studied the case where the network is divided into overlapping subregions, and voltage information of shared buses are duplicated by their connected subregions. In \cite{erseghe2014distributed}, ADMM is directly applied to the underlying distributed reformulation, and convergence is established by assuming nonconvex OPF and ADMM subproblems have zero duality gaps. In \cite{erseghe2015distributed}, standard techniques  {used in the Augmented Lagrangian Method (ALM)} are adopted inside ADMM. Subsequential convergence is proved under the assumption that the penalty parameter stays finite, which is basically assuming the algorithm converges to a feasible solution. The assumption is quite strong as quadratic penalty in general does not admit an exact penalization for nonconvex problems \cite{feizollahi2017exact}. A more recent work \cite{mhanna_adaptive_2019} applied ADMM to a component-based distributed reformulation of AC OPF. The ADMM penalty is adaptively changed in every iteration, and the resulting algorithm numerically converges for various networks under adaptive hyperparameter tuning. In summary, existing ADMM-based algorithms either directly apply ADMM or its variant as a heuristic, or rely on strong assumptions to establish asymptotic convergence. See \cite{molzahn2017survey} for a recent survey on distributed optimization techniques for OPF.

{Another related work is the Augmented Lagrangian Alternating Direction Inexact Newton (ALADIN) algorithm proposed in \cite{houska2016augmented} and \cite{engelmann2018toward}, {which is a mixture of sequential quadratic programming (SQP) and ADMM. The algorithm requires a centralized consensus step that solves an equality-constrained nonconvex quadratic program, which uses the Hessian information of an augmented Lagrangian function. The authors established global convergence of ALADIN to an approximate stationary point when each agent's subproblem is solved to global optimality. However, if subproblems cannot be solved to global optimum due to nonconvexity or numerical considerations, then a good initial point close to a stationary point of the original problem is needed and only local convergence of ALADIN is guaranteed. Due to the use of second-order Hessian information in an SQP framework, ALADIN enjoys local quadratic convergence properties under some technical assumptions.}
}

\subsection{Contribution}
In this paper, we address the convergence issues of ADMM by proposing a two-level distributed algorithmic framework. The proposed framework is motivated by our observation that some crucial structure necessary for the convergence of nonconvex ADMM is absent in traditional distributed reformulations of AC OPF, and we overcome such technical difficulty by embedding a three-block ADMM inside the classic ALM framework. We present the global convergence {to a stationary point} and iteration complexity results of the proposed framework, which rely on mild and realistic assumptions. We demonstrate the convergence, scalability, and robustness of the proposed algorithm over some largest test cases from NESTA \cite{coffrin2014nesta} and PGLib-OPF \cite{babaeinejadsarookolaee2019power}, on which existing ADMM variants may fail to converge. Generically, distributed algorithms can be slow due to limited access to global information and communication delay; however, we show that, with proper parallel implementation, the proposed algorithm achieves fast convergence close to or even better than centralized solver.
    
\subsection{Notation and Organization}
Throughout this paper, we use $\R^n$ to denote the $n$-dimensional real Euclidean space; the inner product of $x,y\in \R^n$ is denoted by $\langle x, y \rangle$; the Euclidean norm is denoted by $\|x\|$, and the $\ell_\infty$ norm is denoted by $\|x\|_\infty$. When $x$ consists of $p$ subvectors, we write $x = (x_1,\cdots, x_p)$. {For a matrix $A\in \R^{m \times n}$, we use $\mathrm{Im}(A)$ to denote its column space}. We use $\Z_{++}$ to denote the set of positive integers, and $[n]=\{1,\cdots, n\}$. For a closed set $C\subset \R^n$, the {orthogonal} projection onto $C$ is denoted by $\mathrm{Proj}_C(x)$, and the indicator function of $C$ is denoted by $\delta_C(x)$, which takes value 0 if $x\in C$ and $+\infty$ otherwise.

The rest of this paper is organized as follows. In section \ref{sec: background}, we review the AC OPF problem and nonconvex ADMM literature. Then in section \ref{sec: reformulation}, we present a new distributed reformulation and the proposed two-level algorithm. In section \ref{sec: convergence}, we state the main convergence results of the proposed two-level algorithm, and discuss related convergence issues. Finally, we present computational experiments in section \ref{sec: numerical}, and conclude this paper in section \ref{sec: conclusion}.

\vspace{-3mm}
\section{Background}\label{sec: background}
\subsection{AC OPF Formulation}
Consider a power network $G(\mathcal{N}, \mathcal{E})$, where $\mathcal{N}$ denotes the set of buses and $\mathcal{E}$ denotes the set of transmission lines. Let $\delta(i)$ be the set of neighbours of $i\in \mathcal{N}$. Let $Y = G+\mathbf{j}B$ denote the complex nodal admittance matrix, where $\mathbf{j} = \sqrt{-1}$ and $G,B\in \R^{|\mathcal{N}|\times|\mathcal{N}|}$. Let $p^g_i$, $q^g_i$ (resp. $p^d_i$, $q^d_i$) be the real and reactive power produced by generator(s) (resp. loads) at bus $i$; if there is no generator (resp. load) attached to bus $i$, then $p^g_i$, $q^g_i$ (resp. $p^d_i$, $q^d_i$) are set to 0. The complex voltage $v_i$ at bus $i$ can be expressed by its real and imaginary parts as $v_i = e_i+\mathbf{j}f_i$. The rectangular formulation of AC OPF is given as 
{\begin{subequations}\label{eq: acropf}
\begin{alignat}{2}
	\min~~ & \sum_{i\in \mathcal{N}}f_i(p^g_i)\label{acropf: obj}\\
	\mathrm{s.t.}~~
	& p^g_i - p^d_i =  G_{ii}(e_i^2 + f_i^2) + \notag\\
	&\sum_{j\in {\delta}(i)} G_{ij} (e_ie_j+f_if_j)-B_{ij}(e_if_j-e_jf_i), &\forall i \in \mathcal{N},\label{acropf: kcl_p}\\
	& q^g_i - q^d_i =  -B_{ii}(e_i^2 + f_i^2) + \notag\\
	&\sum_{j\in {\delta}(i)}-B_{ij} (e_ie_j+f_if_j)-G_{ij}(e_if_j-e_jf_i), &\forall i \in \mathcal{N},\label{acropf: kcl_q}\\
	& p_{ij}^2+q_{ij}^2 \leq \overline{s}_{ij}^2, \quad\quad\quad\quad\quad \quad \quad \quad \forall (i, j) \in \mathcal{E},\label{acropf: apparent_flow_limit} \\
	& \underline{v}_i^2\leq e_i^2+f_i^2 \leq \overline{v}_i^2, \quad\quad\quad \quad \quad \quad~  \forall i \in \mathcal{N},\label{acropf: vm_bounds}\\
	& \underline{p^g_i}\leq p^g_i \leq \overline{p^g_i}, ~~\underline{q^g_i}\leq q^g_i \leq \overline{q^g_i}, \quad\quad \forall i \in \mathcal{N}, \label{acropf: q_bounds}
\end{alignat}
\end{subequations}
}where 
{\begin{subequations} \label{eq: flow_on_line}
\begin{align}
	& p_{ij} = -G_{ij}(e_i^2+f_i^2-e_ie_j-f_if_j)-B_{ij}(e_if_j-e_jf_i),\\
	& q_{ij} = B_{ij}(e_i^2+f_i^2-e_ie_j-f_if_j)-G_{ij}(e_if_j-e_jf_i).
\end{align}
\end{subequations}
}In \eqref{acropf: obj}, the objective $f_i(p^g_i)$ represents the real power generation cost at bus $i$. Constraints \eqref{acropf: kcl_p} and \eqref{acropf: kcl_q} correspond to real and reactive power injection balance at bus $i$.  The real and reactive power flow $p_{ij}, q_{ij}$ on line $(i,j)$ are given in \eqref{eq: flow_on_line}, and \eqref{acropf: apparent_flow_limit} restricts the apparent power flow on each transmission line. Constraints \eqref{acropf: vm_bounds}-\eqref{acropf: q_bounds} limit voltage magnitude, real power output, and reactive power output at each bus to its physical capacity. Since the objective is typically linear or quadratic in real generation, formulation \eqref{eq: acropf} is a nonconvex quadratically constrained quadratic program (QCQP) problem.

\subsection{Nonconvex ADMM}\label{section: ncvx admm}
ADMM was proposed in 1970s \cite{glowinski1975approximation, gabay1976dual} and regarded as a close variant of ALM \cite{hestenes1969multiplier, powell1967method}. The standard ADMM framework consists of a Gauss-Seidel type update on blocks of variables in each ALM subproblem and then a dual update using current primal residuals. The update of each individual block can be decomposed and carried out in parallel given that certain separable structures are available. 
More recently, researchers{\cite{wang2014convergence, wang2015global, jiang2019structured}} realized that the ADMM framework can be used to solve more complicated  nonconvex multi-block problems in the form
\begin{subequations}
	\begin{align}
		\min_{x = (x_1,\cdots, x_p)} \quad & \sum_{i=1}^p f_i(x_i) + g(x)\\
		\mathrm{s.t.}\quad & \sum_{i=1}^p A_ix_i = b,~~ x_i \in \mathcal{X}_i~~\forall i\in[p],
	\end{align}
\end{subequations}
where {there are $p$ blocks of variables $x_i\in \R^{n_i}$ for $i\in [p]$}, and $f_i$'s, $\mathcal{X}_i$'s, and $g$ can be potentially nonconvex. Different assumptions on the problem data are proposed to ensure global convergence to stationary solutions and in general an iteration complexity of $\mathcal{O}(1/\epsilon^2)$ is expected. Though motivated by different applications and adopting different analysis techniques, all these convergence results on nonconvex ADMM rely on the following three assumptions:
\begin{itemize}
	\item[(a)] All nonconvex subproblems need to be solved to global optimality;
	\item[(b)] The functions $f_p$ and $g$ are Lipschitz differentiable, and $\mathcal{X}_p = \R^{n_p}$;
	\item[(c)] The column space of the last block coefficient matrix $A_p$ is sufficiently large, i.e., $\mathrm{Im}([A_1,\cdots, A_{p-1}, b])\subseteq \mathrm{Im}(A_p)$, {where $[A_1,\cdots, A_{p-1}, b]$ is the matrix concatenated by columns of $A_1,\cdots, A_{p-1}$ and $b$.}
\end{itemize}

These three assumptions together restrict the application of ADMM on the OPF problem. No matter what reformulation of \eqref{eq: acropf} is used, ADMM subproblems would still have highly nonconvex constraints, so Assumption (a) is unrealistic for the OPF problem. 
{Assumption (b) is used to provide control of dual variables using primal variables in ADMM. This control is necessary for the construction of a potential function in the convergence analysis. Assumption (c) is needed to guarantee feasibility. If Assumption (c) is not satisfied, then it is possible that ADMM will converge to some $x_1^*,\cdots, x^*_{p-1}$ such that the linear system $A_px_p = b-\sum_{i=1}^{p-1}A_ix_i^*$ has no solution, and hence ADMM fails to find a feasible solution.}
It turns out that Assumptions (b) and (c) cannot be satisfied simultaneously if ADMM were to achieve parallel computation among different agents \cite{sun2019two}. Such limitation motivates us to go beyond the framework of ADMM.

\section{A New Distributed Reformulation and a Two-level ADMM Algorithm}\label{sec: reformulation}

\subsection{A New Distributed Reformulation}
Suppose the network $G$ is partitioned into $R$ subregions $\mathcal{R}_1,\cdots, \mathcal{R}_R\subseteq \mathcal{N}$, each of which is assigned to a local control or operating agent, {i.e., $\mathcal{R}_i$'s are disjoint and $\cup_{i=1}^R \mathcal{R}_i =\mathcal{N}$.} We say $(i,j)\in \mathcal{E}$ is a tie-line if $i\in \mathcal{R}_r$, $j\in \mathcal{R}_l$, and $r\neq l$. Agent $r$ controls variables $x_i = (p^g_i, q^g_i, e_i, f_i)$ for all $i\in \mathcal{R}_r$.  {We say $i\in\mathcal{R}_r$ is a \textit{boundary bus} of $\mathcal{R}_r$ if $i$ is connected to another subregion through a tie-line, and denote the set of boundary buses in $\mathcal{R}_r$ by $ B(\mathcal{R}_r)$.}
We extend the notation $\delta(\mathcal{R}_r)$ to denote the set of all buses connected to (but not in) $\mathcal{R}_r$ by some tie-lines.

The constraints of OPF couple adjacent agents. For example, suppose $(i,j)$ is a tie-line where $i\in \mathcal{R}_r$ and $j\in \mathcal{R}_l$. Agent $r$ requires the information of variables $(e_j, f_j)$ to construct constraints \eqref{acropf: kcl_p}-\eqref{acropf: apparent_flow_limit}; however, agent $r$ cannot directly access $(e_j, f_j)$ as they are controlled by agent $l$, and agent $l$ faces the same situation. In order for these two agents to solve their localized problems in parallel, it is necessary to break the coupling by introducing auxiliary variables. We let each agent $r$ keep additional variables $x_j^r = (e^r_j, f^r_j)$ for all $j\in \delta(\mathcal{R}_r)$. A direct consequence is that all constraints in formulation \eqref{eq: acropf} are decomposed to local agents. For example, for each $i\in \mathcal{R}_r$, constraints \eqref{acropf: kcl_p}-\eqref{acropf: apparent_flow_limit} can be rewritten as 
{\begin{subequations}\label{eq: localconstr}
	\begin{align}
		& p^g_i - p^d_i =  G_{ii}(e_i^2 + f_i^2) + \notag\\
		&\sum_{j\in {\delta}(i)\cap \mathcal{R}_r} G_{ij} (e_ie_j+f_if_j)-B_{ij}(e_if_j-e_jf_i)+ \notag \\
		&\sum_{j\in {\delta}(i)\cap \delta(\mathcal{R}_r)} G_{ij} (e_ie^r_j+f_if^r_j)-B_{ij}(e_if^r_j-e^r_jf_i), \\
		& q^g_i - q^d_i =  -B_{ii}(e_i^2 + f_i^2) + \notag\\
	&\sum_{j\in {\delta}(i)\cap \mathcal{R}_r}-B_{ij} (e_ie_j+f_if_j)-G_{ij}(e_if_j-e_jf_i) + \notag \\
	&\sum_{j\in {\delta}(i)\cap \delta(\mathcal{R}_r)}-B_{ij} (e_ie^r_j+f_if^r_j)-G_{ij}(e^r_if_j-e^r_jf_i), \\
	& p_{ij}^2 + q_{ij}^2 \leq \bar{s}_{ij}^2, \quad ~~~\forall (i,j) \in \mathcal{E}, j\in \mathcal{R}_r, \\
	& {p^r_{ij}}^2 + {q^r_{ij}}^2 \leq \bar{s}_{ij}^2, ~\quad \forall (i,j) \in \mathcal{E}, j\not\in \mathcal{R}_r, 
	\end{align}
\end{subequations}}
where $p_{ij}$, $q_{ij}$ are given in \eqref{eq: flow_on_line}, and
{ \begin{subequations} \label{eq: flow_on_line_r}
\begin{align}
	& p^r_{ij} = -G_{ij}(e_i^2+f_i^2-e_ie^r_j-f_if^r_j)-B_{ij}(e_if^r_j-e^r_jf_i),\\
	& q^r_{ij} = B_{ij}(e_i^2+f_i^2-e_ie^r_j-f_if^r_j)-G_{ij}(e_if^r_j-e^r_jf_i).
\end{align}
\end{subequations}}Notice that all variables appeared in \eqref{eq: localconstr} are controlled by agent $r$, and all such variables are grouped together and denoted by $x^r = \left( \{x_i\}_{i \in \mathcal{R}_r}, \{x^r_j\}_{j\in \delta(\mathcal{R}_r)}\right)$. Moreover, local nonconvex constraints of $\mathcal{R}_r$ can be conveniently expressed as 
\begin{align}\label{eq: localconstr_set}
	\mathcal{X}_r = \{x^r~|~\eqref{acropf: vm_bounds}- \eqref{acropf: q_bounds}, \eqref{eq: localconstr}~\forall i \in \mathcal{R}_r\},
\end{align}
{where, allowing a minor clash of notation, \eqref{acropf: vm_bounds} and \eqref{acropf: q_bounds} are meant to be satisfied for all $i\in \mathcal{R}_r$ in \eqref{eq: localconstr_set}. }
Notice that for every $j \in \cup_{r=1}^R \delta(\mathcal{R}_r)$, bus $j$ is connected to some tie-line, and thus at least two regions need to keep a local copy of $(e_j, f_j)$. We use $R(j)$ to denote the subregion where bus $j$ is located, and $N(j)$ to denote the set of subregions that share a tie-line with $R(j)$ through bus $j$. Naturally we want to impose consensus on local copies of the same variables:
\begin{equation}\label{eq: couple1}
	e^l_j = e_j, ~f^l_j = f_j, ~~ \forall l\in N(j).
\end{equation}
{The regional decoupling techniques used in \eqref{eq: localconstr}-\eqref{eq: couple1} have appeared in the early work by Kim and Baldick \cite{kim1997coarse} among others, where the authors applied a linearized proximal ALM to a distributed OPF formulation. A graphical illustration can be found in \cite{kim1997coarse}.} 

{When ADMM is considered, agents from $N(j)$ and agent $R(j)$ will need to alternatively solve their subproblems in order to parallelize the computation. As we will explain in Section \ref{sec: divergent}, ADMM does not guarantee convergence when both subproblems carry nonconvex functional constraints.} In order to solve this issue, we follow the idea proposed in \cite{sun2019two} by using a global copy $\bar{x}_j = (\bar{e}_j, \bar{f}_j)$ and local slack variables $z^l_j = (z^{l}_{e_j}, z^l_{f_j})$ for $l\in N(j)\cup \{R(j)\}$. For notational consistency, we also write $x^{R(j)}_j = (e_j, f_j) = (e_j^{R(j)}, f_j^{R(j)})$. The consensus is then achieved through
\begin{subequations}\label{eq: couple2}
\begin{alignat}{2}
	e^l_j - \bar{e}_j + z^l_{e_j} =& 0, ~~&&  z^l_{e_j} =0,~~\forall l \in N(j)\cup \{R(j)\}, \\
	f^l_j - \bar{f}_j + z^l_{f_j} =& 0, ~~&& z^l_{f_j} =0, ~~\forall l \in N(j)\cup \{R(j)\}.
\end{alignat}
\end{subequations}
Denote $x = \left(\{x^r\}_{r\in [R]}\right)$,  $\bar{x} =\left (\{\bar{x}_j\}_{j\in\cup_r \delta(\mathcal{R}_r)}\right)$, and $z = \left( \{z^l_j\}_{l\in N(j)\cup\{R(j)\}, j \in \cup_r\delta(\mathcal{R}_r)}\right)$. {Notice that $\bar{x}$ and $z$ are only introduced for boundary buses.} Now we can abstract the AC OPF problem as:
{
\begin{align}\label{eq: distributedOPF}
	\min_{x,\bar{x},z} \quad & \sum_{r=1}^R c_r(x^r):= \sum_{r=1}^R \left(\sum_{i \in \mathcal{R}_r} f_i(p^g_i)\right)\\
	\mathrm{s.t.} \quad & Ax+B\bar{x} +z = 0, \notag \\
	&  x^r \in \mathcal{X}_r~\forall r \in [R], ~\bar{x}\in \bar{\mathcal{X}},  ~z = 0. \notag 
\end{align}
}The objective $c_r(x^r)$ is the sum of all generators' costs in $\mathcal{R}_r$. The linear coupling constraints \eqref{eq: couple2} is compactly expressed as $Ax+B\bar{x}+z = 0$ with matrices $A$ and $B$ of proper dimensions. Each local agent $r$ controls local OPF constraints $\mathcal{X}_r$ defined in \eqref{eq: localconstr_set}. Moreover, without changing the feasible region of \eqref{eq: acropf}, we may restrict $\bar{x}$ inside some convex set $\bar{\mathcal{X}}$. For example, we can simply let $\bar{\mathcal{X}}$ be a hypercube 
\begin{equation}\label{eq: set_bar_x}
	\bar{\mathcal{X}} = \prod_{j\in\cup_r \delta(\mathcal{R}_r)}	 \bar{\mathcal{X}}_j := \prod_{j\in\cup_r \delta(\mathcal{R}_r)}	 \Big\{\bar{x}_j|~\|\bar{x}_j\|_\infty \leq \bar{v}_j\Big\},
\end{equation}
which is compact and easy to project onto.

Next we define stationarity for problem \eqref{eq: distributedOPF}. After projecting out the slack variable $z$, the Lagrangian function of \eqref{eq: distributedOPF} is
\begin{equation}
	L(x,\bar{x}, y) = \sum_{r=1}^R (c_r(x^r)+\delta_{\mathcal{X}_{r}}(x^r)) +\delta_{\bar{\mathcal{X}}}(\bar{x})+\langle y, Ax+B\bar{x}\rangle.
\end{equation}
We use $\partial f(\cdot)$ to denote the general subdifferential of a proper lower semi-continuous function $f:\R^n\rightarrow \R$ \cite[Def 8.3]{rockafellar2009variational}, and $N_{C}(x)$ to denote the general normal cone of $C$ at $x\in C$ \cite[Def 6.3]{rockafellar2009variational}. 
\begin{definition}\label{def: approx_stationary}
	We say $(x,\bar{x}, y)$ is an $\epsilon$-stationary point of problem \eqref{eq: distributedOPF} if there exist $(d_1, d_2, d_3)$ such that $\max\{\|d_1\|, \|d_2\|, \|d_3\|\} \leq \epsilon$ where 
	\begin{subequations}\label{eq: stationary}
		\begin{align}
			d_1 & \in \partial\left(\sum_{r=1}^R c_r(x^r)+\delta_{\mathcal{X}_r}(x^r)\right) +A^\top y, \label{eq: stationary_1} \\
			d_2 & \in N_{\bar{\mathcal{X}}} + B^\top y, \label{eq: stationary_2}\\
			d_3 & = Ax+B\bar{x};
		\end{align}
	\end{subequations}
	or equivalently, $(d_1,d_2,d_3)\in \partial L(x, \bar{x}, y)$. We simply say $(x,\bar{x}, y)$ is a stationary point if $0 \in \partial L(x, \bar{x}, y)$ or $\epsilon=0$.
\end{definition}
If cost functions are concatenated as $c(x) = \sum_{r=1}^R c_r(x^r)$, which is assumed to be continuously differentiable, and let $\mathcal{X} = \prod_{r=1}^R\mathcal{X}_r$, then \eqref{eq: stationary_1} can be further reduced to 
\begin{equation}
	0 \in \nabla c(x) + N_{\mathcal{X}}	(x) + A^\top y.
\end{equation}

\subsection{A Divergent Example for Two-block Nonconvex ADMM}\label{sec: divergent}
\begin{figure}[h!]
		{\centering
		\includegraphics[scale=0.3]{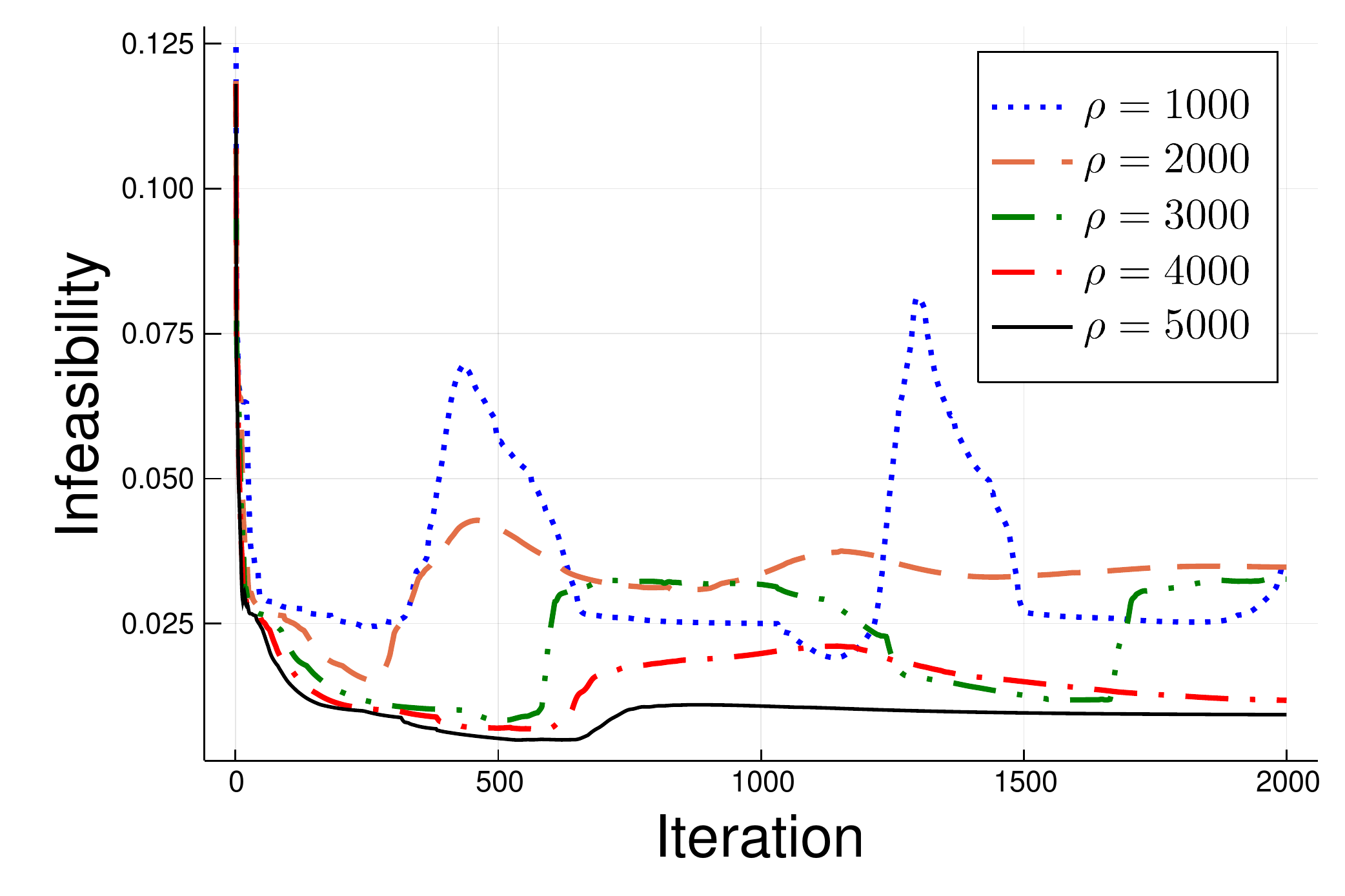}\\
		\includegraphics[scale=0.3]{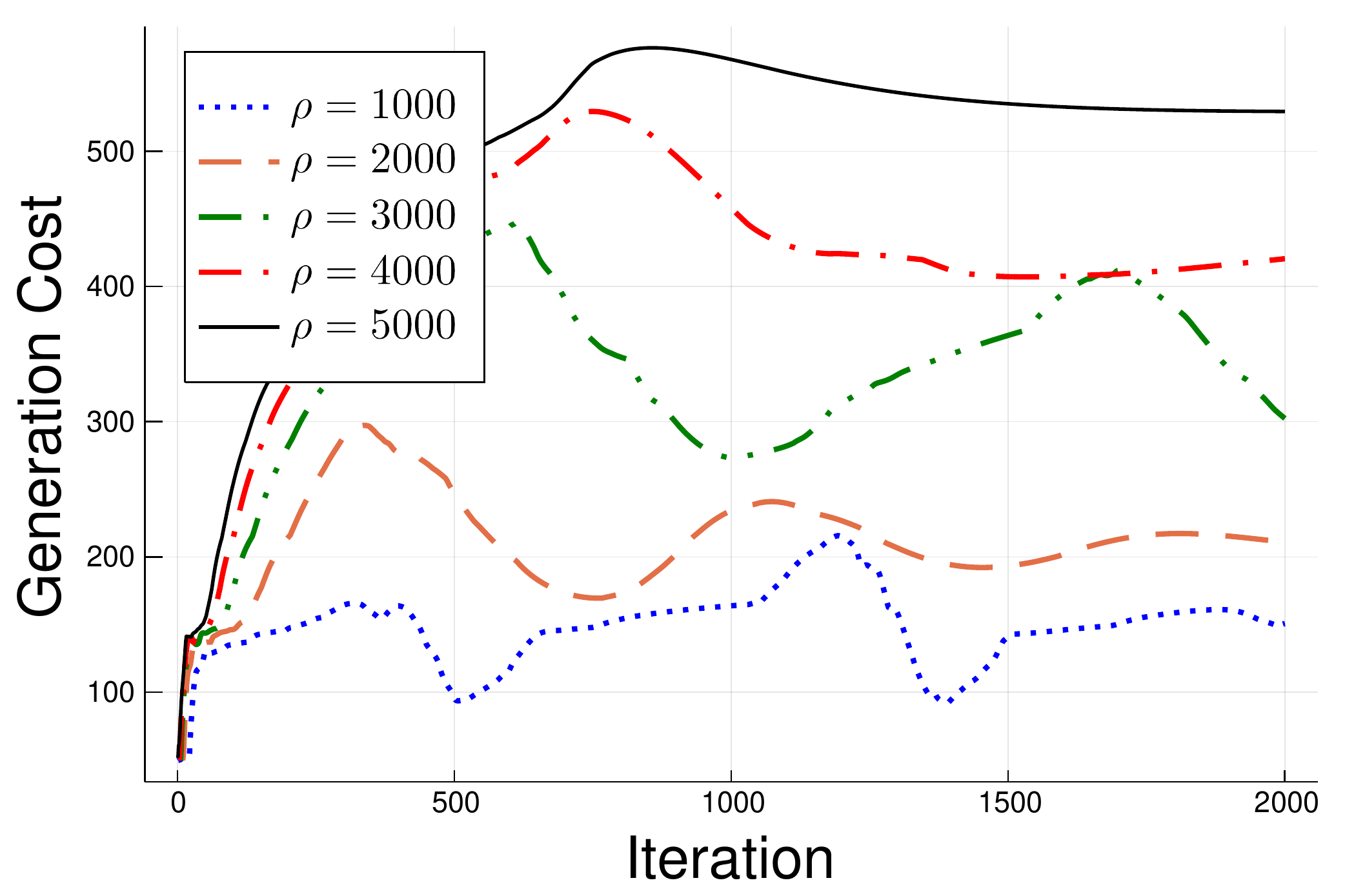}
		\caption{Divergent Behaviour of Vanilla ADMM.}\label{fig:divergent}}
\end{figure}

Before presenting our proposed algorithm, we use a concrete example to demonstrate that the vanilla version of ADMM indeed may suffer divergence for OPF instances. Notice that without introducing slack variable $z$ and the constraint $z=0$ as in \eqref{eq: distributedOPF}, the distributed OPF problem can be formulated as a two-block nonconvex problem 
\begin{align}\label{eq: distributedOPF_2}
	\min_{x \in \mathcal{X},\bar{x}\in \bar{\mathcal{X}}} \quad c(x)\quad  \mathrm{s.t.} \quad  Ax+B\bar{x}  = 0,
\end{align}
which was also directly used to develop distributed algorithms  \cite{erseghe2014distributed, erseghe2015distributed}. We partition the IEEE \texttt{case30} available from \cite{zimmerman2011matpower} into three subregions, directly apply the vanilla version ADMM to the two-block formulation \eqref{eq: distributedOPF_2} with different penalty parameter $\rho$, and plot the \textit{Infeasibility} $\|Ax^t+B\bar{x}^t\|$ and \textit{Generation Cost} $c(x^t)$ as in Figure \ref{fig:divergent}.

As we can see, the primal residual and generation costs exhibit oscillating patterns for $\rho \in \{1000, 2000, 3000, 4000\}$, and do not converge even when ADMM has performed 2000 iterations. For $\rho=5000$, the primal residual converges to 0.0093, and the generation cost at termination is below the lower bound obtained from SOCP; these two facts indicate that ADMM indeed converges to an infeasible solution. 

Such failures of ADMM to obtain feasible solutions result from the tension between the Assumptions (b) and (c) {introduced in Section \ref{section: ncvx admm}}. To be more specific, it is straightforward to verify that $\mathrm{Im}(B) \subset \mathrm{Im}(A)$: given a global copy $\bar{x}_j$, every local agent $l\in N(j)\cup R(j)$ can always keep the same value $x^l_j = \bar{x}_j$ so that the constraint $Ax+B\bar{x}=0$ is satisfied. {Therefore, to satisfy Assumption (c), $A$ need to be the last block and $B$ is the first block in ADMM. But since each agent's problem must consider local OPF constraints, which means $\mathcal{X}$ must be a nonconvex set. This violates the requirement that the problem of the last block must be unconstrained in Assumption (b). Therefore, Assumptions (b) and (c) cannot be satisfied simultaneously. As a result, if we directly apply ADMM to the ACOPF problem, ADMM may fail to converge.} Admittedly, we observe convergence when a even larger penalty $\rho$ is used; however, we want to {emphasize} that despite the numerical success with adaptive parameter tuning \cite{mhanna_adaptive_2019}, the traditional ADMM framework has no guarantee of convergence and could fail for nonconvex distributed OPF.

\subsection{A New Two-level ADMM Algorithm}
In this section we give a full description of the two-level algorithm applied to the OPF problem. The key idea is to dualize and penalize the constraint $z=0$ in \eqref{eq: distributedOPF}, and apply three-block ADMM to solve the {augmented Lagrangian relaxation (ALR)}:
\begin{align}\label{eq: distributedOPF_relaxed}
	\min_{x\in \mathcal{X},\bar{x}\in\bar{\mathcal{X}} ,z} \quad & {c(x)} + \langle \lambda^k, z\rangle + \frac{\beta^k}{2}\|z\|^2 \\
	\mathrm{s.t.} \quad & Ax+B\bar{x} +z = 0, \notag 
\end{align}
with some dual variable $\lambda^k$ and penalty $\beta^k$, to an approximate stationary solution in the following sense.
\begin{definition}\label{def2}
   We say $(x, \bar{x}, z, y)$ is an $\epsilon$-stationary point of problem \eqref{eq: distributedOPF_relaxed} if there exists $(d_1, d_2, d_3)$ such that $\max\{\|d_1\|, \|d_2\|, \|d_3\|\} \leq \epsilon$, where $d_1$ and $d_2$ satisfy \eqref{eq: stationary_1}-\eqref{eq: stationary_2}, $d_3 = Ax+B\bar{x}+z$, and $\lambda^k+\beta^kz+y = 0$.
\end{definition}
Then at termination of ADMM, we update $\lambda^{k+1}$ and $\beta^{k+1}$ as in the classic ALM framework in order to drive $z$ to 0. {To summarize, the proposed algorithm consists of two levels: in the inner level (indexed by $t$), we apply ADMM to solve the ALR problem \eqref{eq: distributedOPF_relaxed}; in the outer level (indexed by $k$), we update the dual information $(\lambda^{k+1},\beta^{k+1})$ using the solution $(x^k,\bar{x}^k, z^k)$ returned by ADMM, which falls into the category of ALM, and then restart the next inner level.}
See Algorithm \ref{alg: twolevel} for a detailed description.

\begin{algorithm}[!htb]
	\caption{: A Two-level ADMM Algorithm}\label{alg: twolevel}
	\begin{algorithmic}[1]
		\STATE \textbf{Initialize} starting points $(x^0, \bar{x}^0)\in \mathcal{X}\times \bar{\mathcal{X}}$ and $\lambda^1\in [\underline{\lambda}, \overline{\lambda}]$; $\beta^1>0$, $k \gets 1$;
		\WHILE{outer stopping criteria is not satisfied}
		\STATE \textbf{Initialize}  $(x^0, \bar{x}^0, z^0, y^0)$ such that $\lambda^k +\beta^kz^0+y^0=0$; $\rho \gets 2\beta^k$, $t\gets 1$;
		\WHILE{inner stopping criteria is not satisfied}\label{line: admm start}
		\STATE each agent $r\in [R]$ updates $(x^r)^{t+1}$ by solving :
		{\small
		\begin{align}\label{eq: subproblem_opf}
			\min_{x^r\in \mathcal{X}_r}~& F_r^t(x^r):=  c_r(x^r)  \notag \\
			& + \sum_{j\in \delta(\mathcal{R}_r)\cup B(\mathcal{R}_r)}\Big(  \langle (y^r_j)^t, x^r_j \rangle  \notag \\
			&  + \frac{\rho}{2}\|x^r_j - (\bar{x}_j)^t+ (z^r_j)^t\|^2 \Big);
		\end{align}
		}
		\STATE each agent $r\in [R]$ sends $(x^r_j)^{t+1}$ to $R(j)$ for $j\in \delta(\mathcal{R}_r)$, and receives $\left((x^l_i)^{t+1},(z^l_i)^t ,(y^l_i)^t\right)$ from every agent $l\in N(i)$ for all $i\in B(\mathcal{R}_r);$
		\STATE each agent $r\in [R]$ updates global copy $\bar{x}^{t+1}_i=$ 
		{\footnotesize
		\begin{align}\label{eq: subproblem_x_bar}
			\mathrm{Proj}_{\bar{\mathcal{X}}_i}\left(\frac{\sum_{l\in N(i)\cup R(i)} \left[(y^l_i)^t+\rho\left((x^l_i)^{t+1}+(z^l_i)^t\right)\right]}{(|N(i)|+1)\rho}\right)
		\end{align}
		}for all $i\in B(\mathcal{R}_r);$
		\STATE each agent $r\in [R]$ sends $(\bar{x}_i)^{t+1}$ to agents in $N(i)$ for $i\in B(\mathcal{R}_r)$, and receives $(\bar{x}_j)^{t+1}$ from agent $R(j)$ for all $j\in \delta(\mathcal{R}_r);$
		\STATE each agent $r\in [R]$ update local slack variable  
		{\small
		\begin{align}\label{eq: subproblem_z}
			(z^r_j)^{t+1} = \frac{-(\lambda^r_j)^k-(y^r_j)^t -\rho\left((x^r_j)^{t+1}-(\bar{x}_j)^{t+1}\right)}{\beta^k + \rho}
		\end{align}
		}and dual variable 
		{\small
		\begin{align}\label{eq: subproblem_y}
			(y^r_j)^{t+1} = (y^r_j)^{t} + \rho \left[(x^r_j)^{t+1} -(\bar{x}_j)^{t+1} + (z^r_j)^{t+1}\right]
		\end{align}}for all $j \in \delta(\mathcal{R}_r)\cup B(\mathcal{R}_r)$;
		\STATE $t \gets t+1$;
		\ENDWHILE \label{line: admm end}
		\STATE denote the solution from the inner loop as $(x^{k}, \bar{x}^{k}, z^{k})$;
		
		\STATE each agent $r\in [R]$ updates outer-level dual variable $	(\lambda^{r}_j)^{k+1}$ for all $j \in \delta(\mathcal{R}_r)\cup B(\mathcal{R}_r)$ and penalty $\beta^{k+1}$; \label{alg: outer dual}
		\STATE $k\gets k+1$;
		\ENDWHILE
	\end{algorithmic}
\end{algorithm}

{The inner level is presented in line \ref{line: admm start}-\ref{line: admm end} of Algorithm \ref{alg: twolevel}, where we apply ADMM to solve ALR \eqref{eq: distributedOPF_relaxed}. To facilitate understanding, define the augmented Lagrangian function associated with \eqref{eq: distributedOPF_relaxed} as 
\begin{align}
  L_{\rho}(x, \bar{x}, z&, y) = c(x)+ \langle \lambda^k, z\rangle + \frac{\beta^k}{2}\|z\|^2   \notag \\
	& + \langle y, Ax+B\bar{x}+z\rangle + \frac{\rho}{2}\|Ax+B\bar{x}+z\|^2,	
\end{align}
where $(\lambda^k,\beta^k)$ are considered as parameters. From a centralized point of view, equations \eqref{eq: subproblem_opf}-\eqref{eq: subproblem_z} are the sequential minimization of $L_{\rho}$ with respect to $x\in \mathcal{X}$, $\bar{x}\in \bar{\mathcal{X}}$, and $z$, respectively, which, together with the update of dual variable $y$ in \eqref{eq: subproblem_y}, constitute a single ADMM iteration. Next we describe ADMM from a local point of view. }All agents simultaneously solve lower-dimensional nonconvex subproblems \eqref{eq: subproblem_opf} by some nonlinear optimization solver. Then each agent $r\in [R]$ sends the current local {estimate} of voltage $(x^r_j)^{t+1}$ to $R(j)$ for all neighboring buses $j\in \delta(\mathcal{R}_r)$. For each bus $j$ connected to a tie-line, we let agent $R(j)$ collect {estimates} $\{(x^l_j)^{t+1}\}_{l\in N(j)}$ and update global copy $(\bar{x}_j)^{t+1}$, though in practice any agents from $N(j)$ can be assigned for this task. Notice that the global copy update \eqref{eq: subproblem_x_bar} involves a projection evaluation, which in general does not admit a closed-form solution; however, if \eqref{eq: set_bar_x} is used for $\bar{\mathcal{X}}$, then $\bar{x}^{t+1}_j$ is exactly the component-wise projection of the argument in \eqref{eq: subproblem_x_bar} onto the box $\bar{\mathcal{X}}_j$. After agent $R(j)$ broadcasts $(\bar{x}_j)^{t+1}$ to agents from $N(j)$, all agents are then able to update the slack variable and dual variables as in \eqref{eq: subproblem_z} and \eqref{eq: subproblem_y}. 

When the inner-level iterates satisfy certain stopping criteria, all agents will update the outer-level dual variable $\lambda^{k+1}$ and penalty $\beta^{k+1}$, as in line \ref{alg: outer dual} of Algorithm \ref{alg: twolevel}, which we will elaborate in the next section. 
	

\section{Convergence and Related Issues}\label{sec: convergence}
In this section, we state the convergence results of the two-level ADMM framework for AC OPF, and discuss related issues. We require the following mild assumptions.
\begin{assumption} \label{assumption1}
\begin{itemize}
	\item[(a)] The objective $c_r(\cdot)$ is continuous differentiable. The functional constraints $\mathcal{X}_r$'s and $\bar{\mathcal{X}}$ are compact, and $\bar{\mathcal{X}}$ is convex.
	\item[(b)] For any $t\in \Z_{++}$, every local agent $r\in[R]$ is able to find a stationary solution $(x^r)^{t+1}$ of subproblem \eqref{eq: subproblem_opf} such that $F^t_r\left((x^r)^{t+1}\right)\leq F^t_r\left((x^r)^{t}\right)$.
\end{itemize}
\end{assumption}
For Assumption \ref{assumption1}(a), the objective function $c_r(\cdot)$ is usually linear or convex quadratic with respect to the argument; without loss of generality, we may assume the set $\mathcal{X}_r$ defined in \eqref{eq: localconstr_set} also enforces bounds on local copies $\{(e^r_j, f^r_j)\}_{j\in \delta(\mathcal{R}_r)}$, and thus $\mathcal{X}_r$ is ensured to be compact. The assumption on $\bar{\mathcal{X}}$ is justified in \eqref{eq: set_bar_x}.
{We note that it is possible to allow $c_r$ to be nonsmooth, e.g., piecewise linear, when its general subdifferential $\partial{c_r}$ is well-defined and bounded over $X_r$.}

Assumption \ref{assumption1}(b) requires the nonconvex subproblem \eqref{eq: subproblem_opf} to be solved to a stationary solution $(x^r)^{t+1}$ that has objective value no worse than that of the previous iterate $(x^r)^{t}$. We believe this assumption is reasonable if the nonlinear solver is warm-started with {the previous solution $(x^r)^{t}$ in iteration $t+1$. For example, the nonlinear solver IPOPT \cite{wachter2006implementation} uses a procedure to accept a trial point if the objective value or constraint violation is reduced during its execution. Since we have a feasible solution $(x^r)^{t}\in \mathcal{X}_r$ to start with, it is reasonable to expect some improvement in the objective. We note that Assumption 1(b) is imposed on the \textit{solution oracle} of subproblem (16), and does not impose any restriction on the initial point $(x^0,\bar{x}^0)\in \mathcal{X} \times \bar{\mathcal{X}}$ supplied to the overall two-level algorithm, which, therefore, still enjoys global convergence as shown below.} 
In addition, Assumption \ref{assumption1}(b) is much weaker than assuming $(x^r)^{t+1}$ is a local or global minimizer of \eqref{eq: subproblem_opf}, which is a common assumption in the literature on nonconvex ADMM. 

\subsection{Global Convergence}
We consider two different rules for updating $(\lambda^{k+1},\beta^{k+1})$. Let $c>1$, $\theta \in [0,1)$, and $\{\eta_k\}_k$ be a nonnegative sequence convergent to 0:
\begin{subequations}\label{eq: condition1}
\begin{align}
    \lambda^{k+1} = & \mathrm{Proj}_{[\underline{\lambda},\overline{\lambda}]}(\lambda^k+\beta^kz^k),\label{eq: outerdual update}\\
    \beta^{k+1} = &
	    	\begin{cases}
		    \beta^k   & \text{if~} \|z^k\|\leq \theta \|z^{k+1}\, \\
		    c\beta^k  &  \text{otherwise,}
		\end{cases}
\end{align}
\end{subequations}
and
\begin{align}\label{eq: condition2}
		(\lambda^{k+1}, \beta^{k+1}) = 
		\begin{cases}
			(\lambda^k+\beta^kz^k, \beta^k) & \text{if~} \|z^k\|\leq \eta_k, \\
			(\lambda^k, c\beta^k)  &  \text{otherwise.}
		\end{cases}
	\end{align}
\begin{theorem}[Global Convergence]\label{thm:globalconv}
	Suppose Assumption \ref{assumption1} holds, and the $k$-th inner-level ADMM is solved to an $\epsilon_k$-stationary point $(x^k, \bar{x}^k, z^k, y^k)$ of \eqref{eq: distributedOPF_relaxed} such that $\epsilon_k\rightarrow 0$ as $k\rightarrow +\infty$. Moreover, the outer-level dual variable $\lambda^{k+1}$ and penalty $\beta^{k+1}$ are updated according to either \eqref{eq: condition1} or \eqref{eq: condition2}. Then the following claims hold.	
	\begin{enumerate}
		\item The sequence $\{(x^k, \bar{x}^k, z^k)\}_k$ is bounded, and therefore there exists at least one limit point $(x^*,\bar{x}^*, z^*)$, where $x^*\in \mathcal{X}=\prod_{r=1}^R \mathcal{X}_r$ and $\bar{x}^*\in \bar{X}$.
		\item Either $(x^*,\bar{x}^*)$ is feasible, i.e., $Ax^*+B\bar{x}^*=0$, or $(x^*,\bar{x}^*)$ is a stationary point of the feasibility problem 
		\begin{align}\label{eq: feasibility}
			\min_{x,\bar{x}}~\frac{1}{2}\|Ax+B\bar{x}\|^2~ \mathrm{s.t.}~x\in \mathcal{X}, \bar{x}\in \bar{\mathcal{X}}.
		\end{align} 
		\item Suppose problem \eqref{eq: distributedOPF} is feasible and the set of stationary points is nonempty. Let $(x^*,\bar{x}^*, z^*)$ be a limit point, and $\{(x^{k_r}, \bar{x}^{k_r}, z^{k_r})\}_r$ be the subsequence convergent to it. If $\{y^{k_r}\}_r$ has a limit point $y^*$, then $(x^*, \bar{x}^*, y^*)$ is a stationary point of problem \eqref{eq: distributedOPF}.
	\end{enumerate}
 \end{theorem}
Proof of Theorem \ref{thm:globalconv} is provided in Appendix. We note that in part 3 of Theorem \ref{thm:globalconv}, we make the assumption that the dual variable $\{y^{k_r}\}_r$ has a limit point. This is a standard sequentially bounded constraint qualification (SBQC) \cite{luo1996exact}.

\subsection{Iteration Complexity}
\begin{theorem}[Iteration Complexity]\label{thm:itercomplexity}
	Suppose Assumption \ref{assumption1} holds, and there exists $0 < \bar{L} < +\infty $ such that 
	{\small 
	\begin{align}\label{thm2_assumption}
		\overline{L} \geq & \sum_{r=1}^R c_r((x^r)^0) + \langle \lambda^k, z^0\rangle + \frac{\beta^k}{2}\|z^0\|^2 \notag \\
		& + \langle y^0, Ax^0+B\bar{x}^0+z^0\rangle +\frac{\rho}{2}\|Ax^0+B\bar{x}^0+z^0\|^2
	\end{align}}for all $k \in \Z_{++}$. Further assume each inner-level ADMM is solved to an $\epsilon$-stationary point $(x^k, \bar{x}^k, z^k, y^k)$ of \eqref{eq: distributedOPF_relaxed}, {the outer-level dual variable $\lambda^{k+1}$ is updated by \eqref{eq: outerdual update}, and $\beta^{k+1} = c^{k+1}\beta^0$} for some $c> 1$, $\beta^0>0$.
	Define
	\begin{itemize}
	    \item $\tau = 2\max\{\|A\|, \|B\|, \frac{1}{4\beta^0}\}$,
	    \item $M = \max_{\lambda \in [\underline{\lambda}, \overline{\lambda}]}\|\lambda\|$,
	    \item $\underline{L} = \min_{x\in \mathcal{X}}c(x)-M^2/\beta^0$,
	    \item $r_{\max} = \max_{x\in\mathcal{X}, \bar{x}\in \bar{\mathcal{X}}} \|Ax+B\bar{x}\|$ (since $\mathcal{X}, \bar{\mathcal{X}}$ compact),
	    \item and for $K\in\Z_{++}$, $$T(K) =\left\lceil \left(\frac{4\beta^0(\overline{L}-\underline{L})\tau^2 c}{c-1}\right)\left(\frac{c^{K}-1}{\epsilon^2}\right)\right\rceil +K.$$ 
	\end{itemize}
	Then Algorithm \ref{alg: twolevel} finds an $\epsilon$-stationary solution of problem \eqref{eq: distributedOPF} in no more than 
	$$K_1 = \left\lceil \log_c\left(\frac{2(\overline{L}-\underline{L} +M r_{\max})}{\beta^0\epsilon^2} \right)\right\rceil$$
	outer ALM iterations and $T(K_1) = \mathcal{O}(1/\epsilon^4)$ inner ADMM iterations. 
	Moreover, if {there exists some $\Lambda >0$ such that} $\|\lambda^k + \beta^k z^k\|\leq \Lambda$ for all outer index $k$, then Algorithm \ref{alg: twolevel} finds an $\epsilon$-stationary solution of problem \eqref{eq: distributedOPF} in no more than 
	$$K_2 = \max\Bigg\{ \left\lceil \log_c\left(\frac{1}{\beta^0 \tau} \right)\right\rceil, \left\lceil \log_c\left(\frac{2(\Lambda + M)}{\beta^0\epsilon} \right)\right\rceil \Bigg\}$$
	outer ALM iterations and $T(K_2)= \mathcal{O}(1/\epsilon^3)$ inner ADMM iterations.
\end{theorem}
The proof of Theorem \ref{thm:itercomplexity} is provided in Appendix. We make some remarks.
\begin{enumerate}
    \item The assumption \eqref{thm2_assumption} can be satisfied trivially, for example, if a feasible solution $(x^0,\bar{x}^0)$ for \eqref{eq: distributedOPF} is known a priori and $(x^0, \bar{x}^0, z^0=0)$ is always used to start inner ADMM. In this case, we can choose $\overline{L} = \max_{x\in \mathcal{X}}c(x)$.
    \item In view of \eqref{eq: couple2}, we can calculate $\|A\|$ and $\|B\|$ directly. Each row of $A$ has exactly one non-zero entry, and each column of $A$ has at most one non-zero entry, so $A^\top A$ is a diagonal matrix with either 0 or 1 on the diagonal, and thus $\|A\|=1$. The number of non-zero entries in each column of $B$ specifies how many subregions keep a local copy of this global variable, so $B^\top B$ is also a diagonal matrix, and we have 
	$${\|B\| = \left(\max_{j\in \cup_{r=1}^R B(\mathcal{R}_r)}|N(j)|+1\right)^{1/2}\leq \sqrt{R}.}$$
	\item The complexity results suggest that a smaller $M$ is preferred, and $M=0$ corresponds to the penalty method. However, we empirically observe that a relatively large range for $\lambda$ usually results in faster convergence, which is also better than the $\mathcal{O}(1/\epsilon^3)$ or $\mathcal{O}(1/\epsilon^4)$ iteration upper bound. We believe such phenomena can be explained by the local convergence properties of ALM.  
\end{enumerate}

\section{Numerical Experiments}\label{sec: numerical}
In this section, we demonstrate the performance of the proposed algorithmic framework. All codes are written in the Julia programming language 1.2.0, and implemented on a Red Hat Enterprise Linux Server 7.6 with 85 Intel 2.10GHz CPUs. All nonlinear constrained problems are modeled using the JuMP optimization package \cite{DunningHuchetteLubin2017} and solved by IPOPT with linear solver MA57 \footnote{{The linear solver MA57 is used for both the centralized algorithm (IPOPT) and the proposed distributed algorithm. It is an interesting research question to fully test IPOPT with parallel linear solvers such as Pardiso or MA97.}}.

\subsection{Network Information and Partition Generation}

{We experiment on four networks: \texttt{case9241\_pegase} (9K), \texttt{case13659\_pegase} (13K) from NESTA \cite{coffrin2014nesta}, and  \texttt{case2848\_rte} (2K), \texttt{case30000\_goc} (30K) from PGLib-OPF \cite{babaeinejadsarookolaee2019power}. See Table \ref{table0} for centralized information.}
\begin{table}[!ht]
	{\footnotesize
		\caption{{Centralized Information.}}\label{table0}
\begin{center}
\begin{tabular}{c|c@{\hskip 0.05in}c@{\hskip 0.05in}|c@{\hskip 0.05in}c@{\hskip 0.05in}|c}
\toprule
	Case	&	AC Obj.	&	AC Time(s)	& SOCP Obj.	& SOCP Time(s)&	Gap (\%) \\
\hline
9K	&	315913.26	&	51.13&	310382.99	&	1224.81 & 1.76\\
13K	&	386117.10	&  160.10&	380262.34	&	990.50  & 1.51\\
\hline
{2K}  &  {1286608.20}  &    {8.75}& {1285025.40}  &    {18.62}  &{0.12} \\
{30K} & {1034405.63}   & {1032.92}&  {1003867.51}  &   {256.20}  &{2.95} \\
\bottomrule
\end{tabular}
\end{center}}
\end{table}
{For networks 2K and 30K, we use linear cost for all generators to enhance numerical stability, which we will elaborate later in Section \ref{sec: implementation}.}
We generate different partitions using the multilevel k-way partitioning algorithm \cite{karypis1998multilevelk} on the underlying graph, which is available from the Julia wrapper of the Metis library \texttt{Metis.jl}. 
{We use the suffix ``-$R$" to indicate that a network is partitioned into $R$ subregions, i.e., 9K-25 refers to network 9K with 25 subregions.}

\subsection{Three Acceleration Heuristics}
It is known that the choice of penalty parameter $\rho$ significantly affects the convergence of ADMM and can potentially accelerate the algorithm \cite{mhanna_adaptive_2019}. Indeed, we observed that for some large cases, the inner-level ADMM may suffer from slow convergence to high accuracy when a constant penalty $\rho$ is used. As a result, given parameters $\theta\in [0,1)$ and $\gamma >1$, we propose three different heuristics to properly accelerate Algorithm \ref{alg: twolevel}.
\begin{enumerate}
    \item Adaptive ADMM penalty {(TL-1)}: the inner-level ADMM penalty is indexed as $\rho^t$ and updated as follows: $\rho^{t+1} = \gamma \rho^{t}$ if $\|Ax^{t+1}+B\bar{x}^{t+1}+z^{t+1}\|>\theta \|Ax^{t}+B\bar{x}^{t}+z^t\|$, and $\rho^{t+1}=\rho^t$ otherwise; in words, we increase the ADMM penalty if the three-block residual $\|Ax^{t+1}+B\bar{x}^{t+1}+z^{t+1}\|$ does not decrease sufficiently.
    \item Different ADMM penalties {(TL-2)}: we assign a different ADMM penalty for each row of the coupling constraint $Ax+B\bar{x}+z=0$, and each penalty is updated according to the first heuristic, where the violation of each single constraint is measured, and the corresponding penalty is adjusted.
    Notice the ALM penalty $\beta$ is a fixed constant for all components of $z$ during the inner ADMM.
    \item Different ALM penalties {(TL-3)}: we assign a different ALM penalty $\beta^t_i$ for each component $z_i$ of the slack variable, and also update it inside ADMM iterations: $\beta_i^{t+1} = \gamma \beta_i^{t}$ if $|z_i^{t+1}|>\theta |z_i^t|$, and $\beta_i^{t+1}=\beta_i^t$ otherwise; the corresponding ADMM penalty $\rho^t_i$ is always assigned to be $2\beta^t_i$, as required in our analysis. When the $k$-th ADMM terminates, current values of $z_i^k$ and $\beta_i^k$ are used to update outer level dual variable $\lambda_i^{k+1}$.
\end{enumerate}
We note that the first two heuristics have been used to accelerate ADMM, while the last heuristic also penalizes the slack variable $z$ adaptively in ADMM iterations.
\begin{table*}[t]
 \caption{{Performance of the Two-level ADMM Algorithm on 9K and 13K Networks.}} \label{table1}
\begin{tabularx}{\textwidth}{c@{\hskip 0.05in}c@{\hskip 0.1in}c@{\hskip 0.1in}|c@{\hskip 0.1in}c@{\hskip 0.1in}c@{\hskip 0.06in}c@{\hskip 0.1in}|c@{\hskip 0.1in}c@{\hskip 0.1in}c@{\hskip 0.1in}c@{\hskip 0.1in}|c@{\hskip 0.1in}c@{\hskip 0.1in}c@{\hskip 0.1in}c@{\hskip 0.02in}|c}
\Xhline{2\arrayrulewidth}
      &  & &  \multicolumn{4}{c|}{TL-1} &  \multicolumn{4}{c|}{TL-2} &  \multicolumn{4}{c|}{TL-3}\\
\hline
Case  &  Tie-line & Dim& Outer & Inner  & Gap (\%) & $\|r\|_\infty$  & Outer & Inner  & Gap (\%) & $\|r\|_\infty$ & Outer & Inner  & Gap (\%) & $\|r\|_\infty$ & Avg. Time (s)\\ 
\hline
9K-25&	357 & 2084 & 88 & 317 &	1.96& 2.97E-3&	237& 270&	1.97&	2.89E-3&	246& 261&	1.95&	2.91E-3 & 735.50\\
9K-30&	412 & 2478 & 58 & 189 &	1.51& 3.20E-3&	139& 173&	1.65&	3.49E-3&	148& 163&	1.62&	3.50E-3 & 428.55\\ 
9K-35&	518 & 2984 & 46 & 112 &	1.00& 2.64E-3&	 82& 115&	1.07&	2.63E-3&	 87& 101&	1.10&	2.72E-3 & 158.73\\
9K-40&	514 & 3108 & 57 & 186 &	0.54& 3.91E-3&	118& 151&	0.58&	3.82E-3&	126& 140&	0.57&	3.69E-3 & 271.43\\
9K-45&	603 & 3538 & 39 &  89 &	0.17& 2.13E-3&	 62& 91	&   0.34&	2.18E-3&	 66&  80&	0.28&	2.23E-3 & 121.89\\
9K-50&	676	& 3808 & 40 &  92 &-0.18& 2.57E-3&	 63& 90	&  -0.11&	2.62E-3&	 67&  81&  -0.11&	2.62E-3 & 115.30\\
9K-55&	651	& 3776 & 59 & 181 &	0.47& 4.38E-3&	127& 161&	0.49&	4.28E-3&	133& 148&	0.47&	4.19E-3 & 186.84\\
9K-60&	693	& 4080 & 49 & 120 &	0.00& 3.73E-3&	 85& 118&  -0.02&	3.75E-3&	 90& 104&  -0.06&	3.63E-3 & 126.44\\
9K-65&	741	& 4292 & 55 & 137 &	0.07& 2.94E-3&	106& 140&	0.20&	3.11E-3&	112& 126&	0.20&	3.10E-3 & 152.75\\
9K-70&	764	& 4430 & 38 &  86 &-0.33& 2.26E-3&	 63&  87&  -0.24&	2.34E-3&	 65&  80&  -0.30&	2.29E-3 & 97.04\\
\hline
13K-25	&   371	&	2234	&	32	&	69	&	1.52	&	3.51E-03	&	47	&	73	&	1.58	&	3.38E-03	&	49	&	63	&	1.55	&	3.43E-03 & 1513.66\\
13K-30	&	452	&	2536	&	28	&	57	&	1.30	&	2.26E-03	&	38	&	60	&	1.33	&	2.36E-03	&	39	&	53	&	1.33	&	2.45E-03 & 682.30\\
13K-35	&	482	&	2846	&	32	&	67	&	0.74	&	3.37E-03	&	45	&	67	&	0.77	&	3.45E-03	&	48	&	62	&	0.78	&	3.44E-03 & 887.56\\
13K-40	&	533	&	3066	&	25	&	48	&	1.20	&	2.73E-03	&	34	&	54	&	1.30	&	2.73E-03	&	35	&	48	&	1.28	&	2.69E-03 & 533.48\\
13K-45	&	655	&	3768	&	31	&	64	&	1.19	&	4.05E-03	&	43	&	66	&	1.28	&	4.07E-03	&	46	&	60	&	1.27	&	4.05E-03 & 300.73\\
13K-50	&	618	&	3692	&	25	&	48	&	0.53	&	2.82E-03	&	33	&	56	&	0.60	&	3.03E-03	&	35	&	49	&	0.59	&	3.15E-03 & 330.53\\
13K-55	&	721	&	4246	&	25	&	46	&	0.62	&	3.98E-03	&	31	&	53	&	0.64	&	4.01E-03	&	33	&	47	&	0.65	&	4.00E-03 & 239.50\\
13K-60	&	717	&	4176	&	22	&	40	&	1.04	&	2.42E-03	&	27	&	47	&	1.07	&	2.46E-03	&	27	&	41	&	1.07	&	2.44E-03 & 157.78\\
13K-65	&	736	&	4258	&	21	&	38	&	0.81	&	2.11E-03	&	26	&	44	&	0.82	&	2.18E-03	&	26	&	40	&	0.83	&	2.23E-03 & 197.60\\
13K-70	&	843	&	4784	&	25	&	47	&	1.19	&	4.24E-03	&	32	&	56	&	1.25	&	4.16E-03	&	29	&	53	&	1.22	&	4.21E-03 & 167.57\\
\Xhline{2\arrayrulewidth}
\end{tabularx}
\end{table*}
\subsection{Implementation Details}\label{sec: implementation}
\subsubsection*{Parallelization of nonconvex subproblems} Each JuMP  model carrying a subregion's localized OPF constraints is initialized on a {core}. During each (inner) iteration, these models are solved in parallel {on different cores} by IPOPT, which consist of the major computation of the algorithm. {Multiple models on the same core are solved sequentially.}
Then current local solutions are gathered through the master {node}, and auxiliary primal variables and dual variables in different subregions are updated in closed form.

\subsubsection*{Parameters and Initialization} For the heuristics introduced in the previous subsection, we set $\gamma = 6.0$, $\theta = 0.8$; for the first two heuristics, when ADMM terminates, the outer-level penalty is updated as $\beta^{k+1} = c\beta^k$ where $c=6.0$. The initial value $\beta^0$ is set to 1000.0, and an upper bound of $1.0e24$ is imposed in all penalty updates. Each component of $\lambda$ is bounded between $\pm 1.0e12$. Flat start is used to initialize IPOPT: we choose $(e_i, f_i, p^g_i, q_i^g) = (1, 0, 0,0)$ for all $i\in \mathcal{N}$. Dual variables $y^0$ and $\lambda^0$ are initialized with zeros.

\subsubsection*{Scaling of IPOPT} The proposed algorithm inevitably needs to deal with potentially large penalties and dual variables, and we observe that IPOPT will encounter numerical failures or produce oscillating solutions. To overcome this problem, we manually scale the objective of the each JuMP model so that the largest coefficient passed to the solver is in the order of $1.0e+8$. This trick helps IPOPT output stable solutions efficiently. 
{Nevertheless, we observe that the proposed algorithm tends to yield large gaps on instances with quadratic generation cost. Since we scale the objective inside the execution of IPOPT, the numerical inaccuracy will be magnified when we calculate the true generation cost; in this situation, a quadratic cost function is more sensitive than a linear one. Such numerical issues are known to be associated with penalty-type methods and deserve further investigations.}

\subsubsection*{Termination of Inner and Outer Iterations} We stop the inner-level ADMM if: (1) $\|Ax^t+B\bar{x}^t+z^t\|\leq \sqrt{d}/(2500k)$ where $d$ is the number of the coupling constraints and $k$ is the current outer-level iteration index, \textbf{or} (2) $\|z^{t}-z^{t-1}\|\leq 1.0e-8$. The first condition ensures that the ADMM primal residual is under certain tolerance, which also tends to 0 as $k\rightarrow \infty$. The second condition measures the dual residual of the last block in \eqref{eq: distributedOPF_relaxed}; if this quantity is small, then we believe $z^t$ has stabilized in the current inner iteration, which encourages us to terminate ADMM early and proceed to update outer-level dual variables. {The outer level is terminated if the consensus mismatch satisfies $\|Ax^k+B\bar{x}^k\|\leq \sqrt{d} \times \epsilon$, where $\epsilon >0$ is some given tolerance.}

\subsection{Numerical Performance {on 9K and 13K from NESTA}}
\begin{figure}[h!]
		\centering
		\includegraphics[scale=0.3]{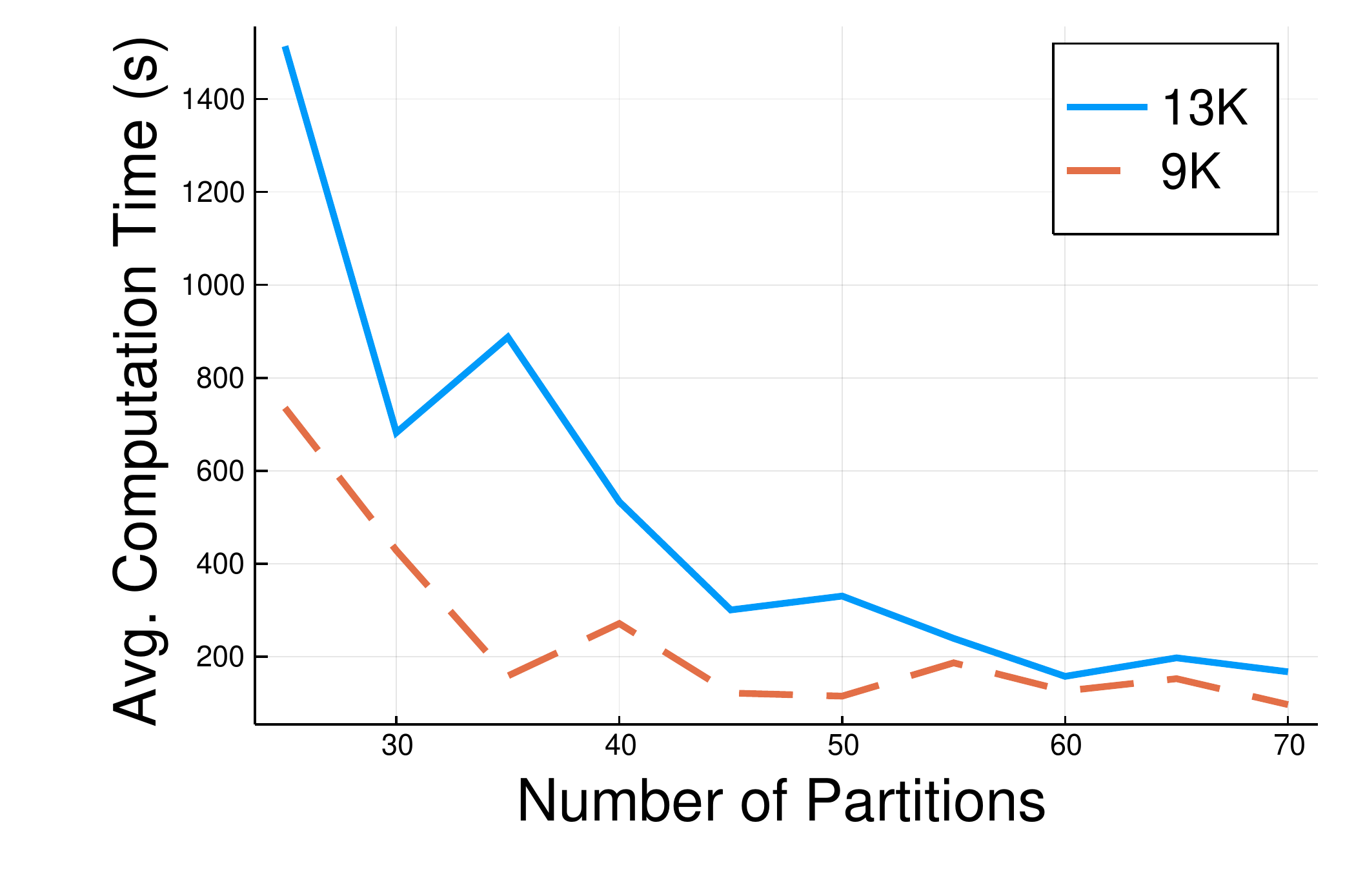}
		\caption{Average Computation Times for 9K and 13K Networks.}\label{fig:time1}
\end{figure}

{The results for 9K and 13K networks are displayed in Table \ref{table1}, where we set $\epsilon$ = $2\times10^{-4}$. We let each core solve a single subproblem (corresponding to a subregion) in every ADMM iteration.} The total numbers of tie-lines and coupling constraints in each instance are recorded in the second and third columns. For each heuristic, we report the number of outer iterations (Outer), number of inner iterations (Inner), duality gap with SOCP lower bound (Gap (\%)), and max violation of the coupling at termination ($\|r\|_{\infty}$); {the average wall clock time of the three heuristics (Avg. Time (s)) is given in the last column}. The three algorithms reach the desired infeasibility tolerance in all test instances, and the max constraint violation is in the order of $10^{-3}$. 
Overall the generation costs at termination are very close to the SOCP lower bound, indicating the algorithms converge to solutions with high quality. Moreover, we emphasize the scalability of the algorithm by pointing out that, the number of inner and outer iterations are stable across instances, even though the dimension of the coupling ranges from 2000 to near 5000. We plot the averaged computation time for the proposed algorithm combined with three heuristics in Fig. \ref{fig:time1}. The computation time drops significantly and becomes comparable to the centralized solver as the number of partitions increases: 735.50 to 97.04 seconds for 9K-bus system, and 1513.66 to 167.57 seconds for 13K-bus system. 
Our results validate the feasibility of using parallelization to speed up computation.


\begin{table*}[t]
 \caption{{Performance of the Two-level ADMM Algorithm on 2K and 30K Networks.}} \label{table_pglib}
\begin{tabularx}{\textwidth}{c@{\hskip 0.05in}c@{\hskip 0.05in}c|c@{\hskip 0.1in}c@{\hskip 0.1in}c@{\hskip 0.1in}c@{\hskip 0.1in}|c@{\hskip 0.1in}c@{\hskip 0.1in}c@{\hskip 0.1in}c@{\hskip 0.1in}|c@{\hskip 0.1in}c@{\hskip 0.1in}c@{\hskip 0.1in}c@{\hskip 0.05in}|c}
\Xhline{2\arrayrulewidth}
      &  & &  \multicolumn{4}{c|}{TL-1} &  \multicolumn{4}{c|}{TL-2} &  \multicolumn{4}{c|}{TL-3} \\
\hline
Case  &  Tie-line & Dim& Outer & Inner  & Gap (\%) & $\|r\|_\infty$ & Outer & Inner  & Gap (\%) & $\|r\|_\infty$ & Outer & Inner  & Gap (\%) & $\|r\|_\infty$& Avg.Time (s) \\ 
\hline
2K-120  & 693  & 3412  & 124 & 548 & 10.73 & 6.18E-3 & 197 & 573 & 10.75 & 6.14E-3 & 194 & 566 & 10.76 & 6.18E-3 &281.46\\
2K-180  & 857  & 4220  & 140 & 635 & 12.84 & 6.03E-3 & 229 & 671 & 12.91 & 6.03E-3 & 225 & 656 & 12.92 & 6.03E-3 &368.28\\
2K-240  & 1059 & 5190  & 103 & 436 & 3.59  & 7.30E-3 & 141 & 516 & 2.62  & 7.31E-3 & 136 & 512 & 3.10  & 7.32E-3 &291.46\\
2K-300  & 1591 & 7388  & 81  & 308 & -5.94 & 4.78E-3 & 106 & 373 & -5.90 & 4.81E-3 & 101 & 368 & -5.83 & 4.77E-3 &275.90\\
2K-360  & 1345 & 6496  & 93  & 385 & 3.72  & 6.61E-3 & 126 & 456 & 3.32  & 6.62E-3 & 122 & 456 & 3.79  & 6.60E-3 &299.99\\
\hline
30K-2400 & 12761 & 72474 & 106 & 288 & 0.61  & 2.28E-3 & 124 & 363 & 1.15  & 2.27E-3 & 121 & 360 & 1.17  & 2.28E-3 & 2998.99 \\
30K-2700 & 11288 & 66112 & 186 & 633 & 6.18  & 2.23E-3 & 230 & 812 & 6.29  & 2.24E-3 & 226 & 808 & 6.30  & 2.23E-3 & 6937.62 \\
30K-3000 & 10748 & 64768 & 201 & 698 & 6.18  & 1.94E-3 & 250 & 896 & 6.36  & 1.94E-3 & 246 & 893 & 6.38  & 1.93E-3 & 7567.32 \\
30K-3300 & 11802 & 69542 & 151 & 492 & 3.47  & 2.09E-3 & 186 & 631 & 3.60  & 2.06E-3 & 183 & 627 & 3.62  & 2.09E-3 & 5645.78 \\
30K-3600 & 13695 & 78336 & 93  & 255 & 1.58  & 2.11E-3 & 113 & 312 & 1.48  & 2.10E-3 & 110 & 309 & 1.52  & 2.11E-3 & 2946.14 \\
30K-3900 & 14042 & 80332 & 84  & 226 & -0.39 & 2.39E-3 & 103 & 275 & -0.19 & 2.39E-3 & 100 & 275 & -0.18 & 2.39E-3 & 2694.60 \\
30K-4200 & 13415 & 77420 & 122 & 365 & 1.82  & 2.88E-3 & 148 & 468 & 1.84  & 2.87E-3 & 144 & 462 & 1.85  & 2.88E-3 & 4489.54 \\
\Xhline{2\arrayrulewidth}
\end{tabularx}
\end{table*}
\subsection{{Numerical Performance on 2K and 30K from PGLib-OPF}}\label{sec: pglib}
\begin{figure}[h!]
     \centering
     \begin{subfigure}[b]{0.3\textwidth}
         \centering
         \includegraphics[width=\textwidth]{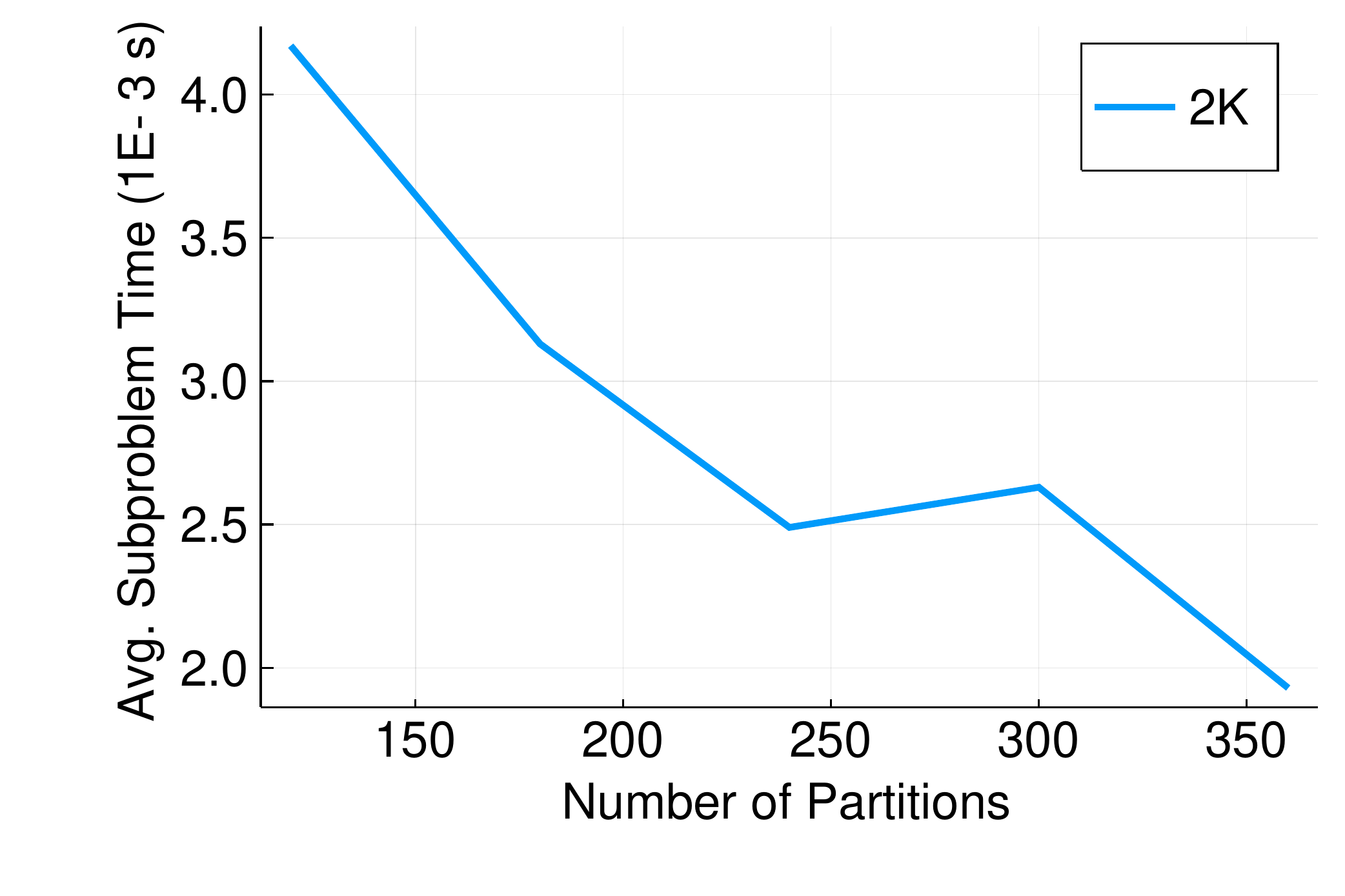}
     \end{subfigure}
     \begin{subfigure}[b]{0.325\textwidth}
         \centering
         \includegraphics[width=\textwidth]{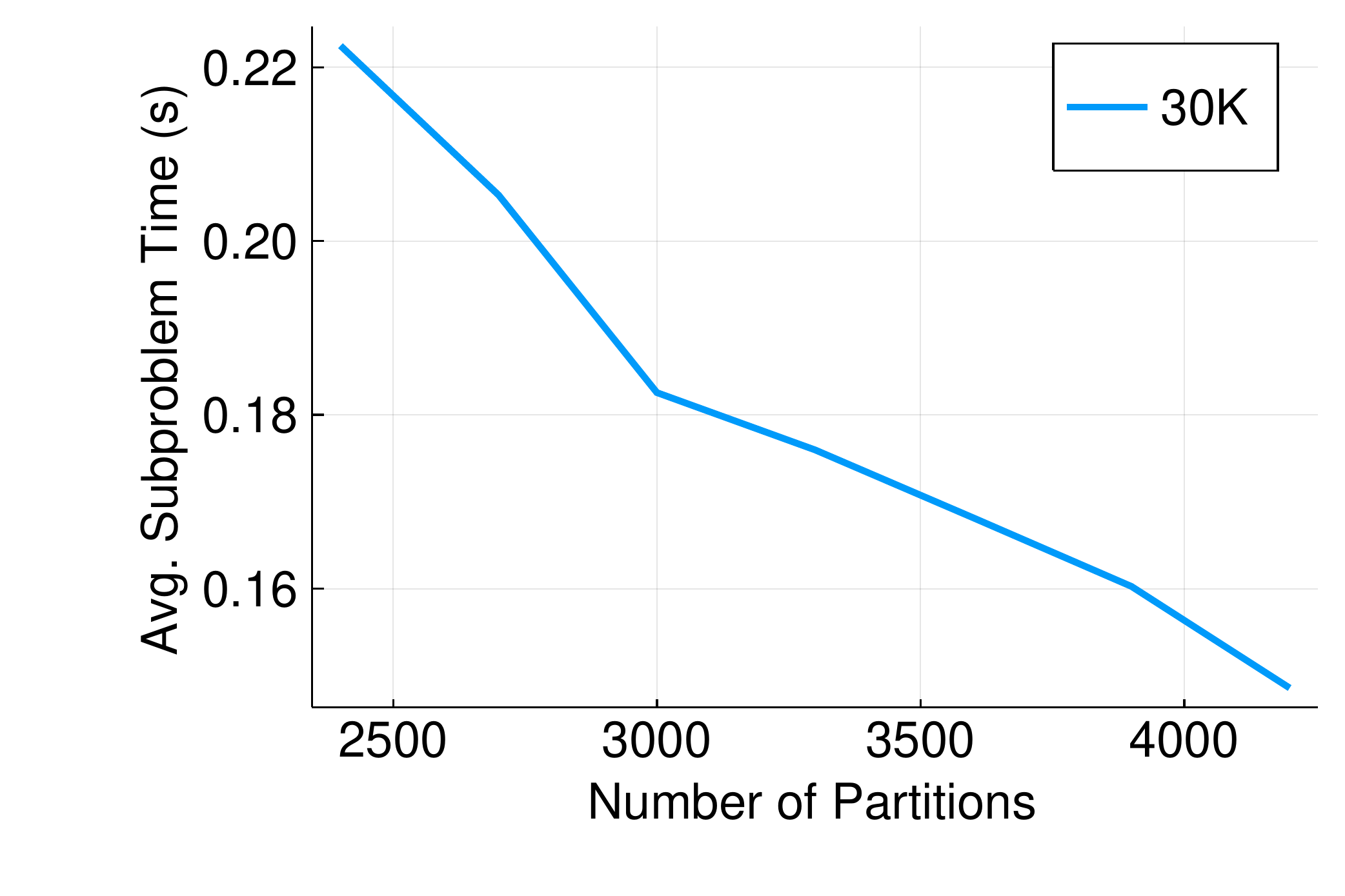}
      \end{subfigure}
      \caption{{Average Subproblem Times for 2K and 30K Networks.}}\label{fig: avg subproblem time}
\end{figure}

{We present the same metrics for networks 2K ($\epsilon$ = $5\times10^{-4}$) and 30K ($\epsilon$ = $10^{-4}$) in Table \ref{table_pglib}. 
We limit the computational resource to 60 cores. Different from experiments in Table \ref{table1}, when the network is partitioned into $R$ subregions, each core solves $R/60$ zonal subproblems sequentially in every inner ADMM iteration. So the computation time presented in Table \ref{table_pglib} is expected to further reduce when more computation power is available. To this end, we plot the averaged subproblem time in Figure \ref{fig: avg subproblem time}, i.e., average inner iteration time$\times 60/R$.}

{We observe that in both Table \ref{table1} and Table \ref{table_pglib}, the duality gap is negative for a few cases. This is because the infeasibility of coupling constrains in $\ell_2$ norm is still relatively large, even though the max violation is already in the order of $10^{-3}$.  As an illustration, we plot the $\ell_2$ and $\ell_\infty$ norm of the infeasibility $r = Ax+B\bar{x}$ over time for 2K-300 in Figure \ref{fig:res}.
\begin{figure}[h!]
     \centering
     \begin{subfigure}[b]{0.31\textwidth}
         \centering
         \includegraphics[width=\textwidth]{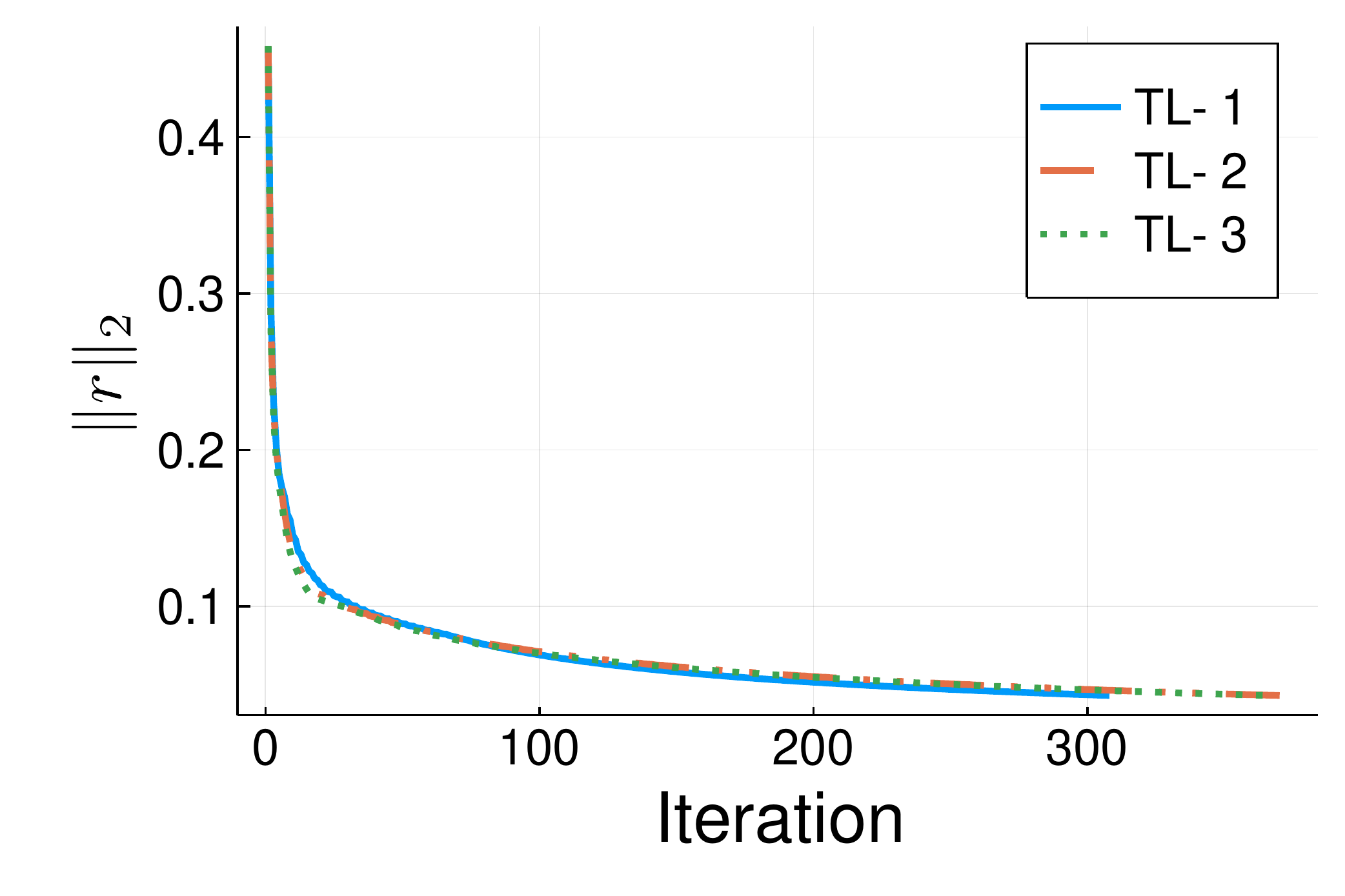}
      \end{subfigure}
     \begin{subfigure}[b]{0.3\textwidth}
         \centering
         \includegraphics[width=\textwidth]{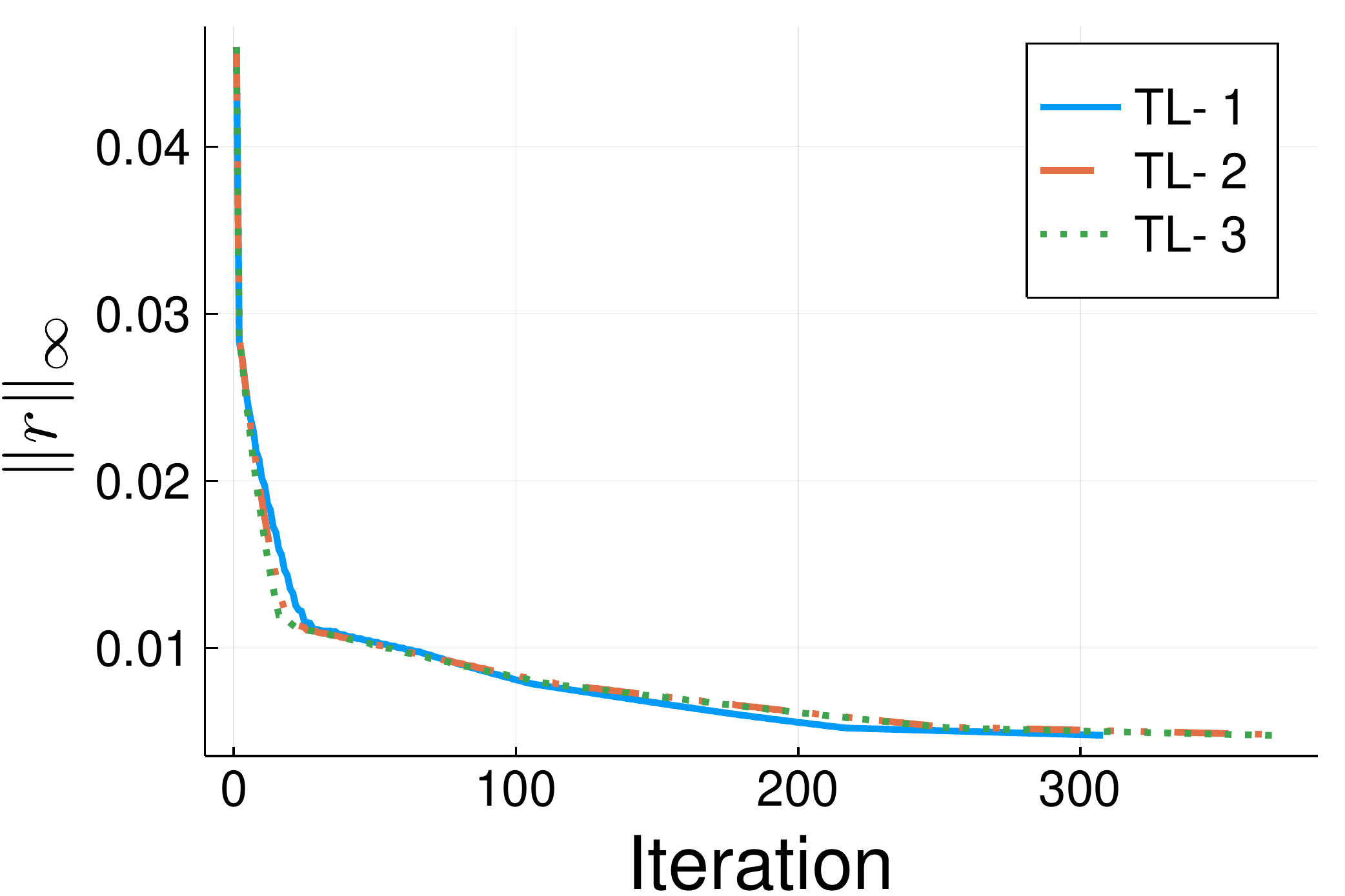}
     \end{subfigure}
        \caption{{Evolution of Primal Residual for 2K-300 Network.}}
        \label{fig:res}
\end{figure}
We note that the numbers of coupling constraints reported in column ``Dim" of Table \ref{table_pglib} are indeed very large, and this is why we choose a smaller tolerance $\epsilon$ = $10^{-4}$ for 30K. This phenomenon suggests that in practice one may need to further decrease $\epsilon$, or let the algorithm run longer to refine local solutions when the network is heavily partitioned.}

\subsection{Comparison with One-level ADMM Variants}
To further illustrate the advantage of the proposed two-level ADMM algorithm, we implement two state-of-the-art one-level ADMM variants and display their results in Table \ref{table_onelevel}. All implementation issues mentioned above are considered; however, for all cases, both one-level ADMM algorithms fail to achieve the desired infeasibility tolerance within 1000 iterations.
For the modified ADMM \cite{erseghe2015distributed}, the penalty parameter quickly reaches the upper bound $10^{24}$, which indicates their assumption for convergence is not satisfied,
and such large penalty leads to solutions with high generation costs. Jiang et al. \cite{jiang2019structured} realize that ADMM can be used to solve the relaxed problem \eqref{eq: distributedOPF_relaxed} with $\lambda^k=0$ and large enough constant $\beta$. When $\beta=\mathcal{O}(1/\epsilon^2)$, their proposed ADMM variant is guaranteed to find an $\epsilon$-stationary solution of the original problem. Using parameters suggested in their analysis, we observe that the true infeasibility  $\|Ax^t+B\bar{x}^t\|$ decreases very slowly, though the three-block residual $\|Ax^t+B\bar{x}^t+z^t\|$ converges to near 0 fast. This indicates the need for a more careful tuning of the penalty $\beta$, but there is no principled way to select all hyperparameters in their algorithm for OPF instances. 
\begin{table}[tbhp]
	\caption{Comparison with One-level ADMM Variants.}\label{table_onelevel}
\begin{center}
\begin{tabular}{c@{\hskip 0.02in}|c@{\hskip 0.02in}c@{\hskip 0.02in}|c@{\hskip 0.02in}c@{\hskip 0.02in}|c@{\hskip 0.02in}c}
\hline 
        & \multicolumn{2}{c}{Modified ADMM\cite{erseghe2015distributed}} & \multicolumn{2}{|c|}{ADMM-g \cite{jiang2019structured}} &
        \multicolumn{2}{c}{Proposed TL-1}\\
\hline
	Case	&	Gap(\%)	&	Time (s)	&	Gap(\%)	&	Time (s) 	&	Gap(\%)	&	Time (s) \\
\hline
13K-25	&	88.98	&	2869.48	&	-6.14	&	1946.21 & 1.52 & 1788.62\\
13K-30	&	90.64	&	2507.26	&	-6.21	&	1874.13 & 1.30 & 680.16\\
13K-35	&	91.86	&	2034.09	&	-6.97	&	1510.31 & 0.74 & 752.41\\
13K-40	&	92.79	&	2168.87	&	-8.09	&	1555.16 & 1.20 & 430.68\\
13K-45	&	93.54	&	1428.20	&	-8.06	&	1436.43 & 1.19 & 306.79\\
13K-50	&	94.15	&	1463.57	&	-8.44	&	1405.47 & 0.53 & 374.15\\
13K-55	&	94.65	&	1521.64	&	-12.59	&	1370.72 & 0.62 & 204.50\\
13K-60	&	95.09	&	1459.83	&	-8.26	&	1208.48 & 1.04 & 157.32\\
13K-65	&	95.44	&	1910.87	&	-9.16	&	1279.62 & 0.81 & 207.58\\
13K-70	&	95.75	&	1821.54	&	-12.20	&	1340.73 & 1.19 & 160.20\\
\hline
\end{tabular}
\end{center}
\end{table}

Finally, we note that the proposed two-level algorithm is not intended to replace centralized solvers or other distributed algorithms, but rather aims to serve as a general algorithmic framework with theoretically supported convergence guarantees, upon which other methods could be applied to further accelerate convergence of distributed computation in various practical applications.

\section{Conclusion}\label{sec: conclusion}
In this paper we propose a two-level distributed ADMM framework for solving AC OPF. The proposed algorithm is proven to have guaranteed global convergence with an iteration complexity upper bound, and can be further accelerated by suitable heuristics. Promising numerical results over some large-scale test cases show that the proposed algorithm provides a new, robust, distributed, and convergence-guaranteed algorithmic framework for solving real-world sized AC OPF problems. {Future directions include extending the two-level ADMM framework to security-constrained OPF problems, and combining with other methods, such as the component-based decomposition and adaptive penalty parameter scheme in \cite{mhanna_adaptive_2019}. Also further experimenting with massive parallelization and studying the impact of communication would be very interesting.}
\\

\noindent\textbf{Acknowledgement}
The authors are grateful for support from the National Science Foundation [Grant ECCS-1751747] and the ARPA-E [Grant DE-AR0001089].

\ifCLASSOPTIONcaptionsoff
  \newpage
\fi

\bibliographystyle{IEEEtran}
\bibliography{ieee_ref}

\section*{Appendix}\label{appendix}
\subsection{Proof of Theorem \ref{thm:globalconv}}
\begin{proof}
The claims under condition \eqref{eq: condition1} are proved in Theorem 1-2 of \cite{sun2019two}. We shall prove the claims under condition \eqref{eq: condition2}.
\begin{enumerate}
\item The first claim follows from the compactness of $\mathcal{X}_r$'s and $\bar{\mathcal{X}}$ and the fact that $\|Ax^k+B\bar{x}^k+z^k\|\rightarrow 0$.
\item 
First we assume the sequence $\{\beta^k\}_k$ stays finite. Then there exists $K>0$ such that $\eta_k \geq \|z^k\|$ for all $k\geq K$, which follows $\|z^k\|\rightarrow 0$. Next assume $\beta^k\rightarrow +\infty$, and assume without loss of generality that $(x^k, \bar{x}^k, z^k)\rightarrow (x^*, \bar{x}^*, z^*)$. If the first case in \eqref{eq: condition2} is executed infinitely many times, we have $\|z^k\|\rightarrow 0$; otherwise we must have $\|\lambda^k\|$ stays constant for all sufficiently large $k$. Define $\tilde{\lambda}^k := \lambda^k+\beta^kz^k=-y^k$. If $\{\tilde{\lambda}^k\}$ has a bounded subsequence, then we have $\|z^k\|\rightarrow 0$ as $\beta^k$ converges to infinity.  
In all previous cases, we have $\|Ax^k+B\bar{x}^k\| \leq \|Ax^k+B\bar{x}^k+z^k\|+\|z^k\|\leq \epsilon_k + \|z^k\|\rightarrow 0$, and thus $\|Ax^*+B\bar{x}^*\|=0$. In the last case we have $\beta^k\rightarrow\infty$, $\|\lambda^k\|$ stays bounded, and $\|\tilde{\lambda}^k\|\rightarrow \infty$. By the definition of $\tilde{\lambda}^k$, we can see $y^k/\beta^k \rightarrow -z^*=Ax^*+B\bar{x}^*$. At termination of ADMM, we have
	{
	\begin{subequations}\label{eq: admm_stationary}
	\begin{align}
		d_1^k &\in \nabla c(x^k) + A^\top y^k + N_{\mathcal{X}}(x^k)\label{eq: outer_stat_1}\\
		d_2^k &\in B^\top y^k +N_{\bar{\mathcal{X}}}(\bar{x}^k)	\label{eq: outer_stat_2}\\
		d_3^k &= Ax^k+B\bar{x}^k+z^k, 
	\end{align}
	\end{subequations}}where $\max \{\|d_1^k\|, \|d_2^k\|, \|d_3^k\|\}\leq \epsilon_k\rightarrow0$. By the closeness of normal cone, dividing vectors in \eqref{eq: outer_stat_1}-\eqref{eq: outer_stat_2} by $\beta^k$, and then taking limit, we have $0 \in  A^\top(Ax^*+B\bar{x}^*) + N_{\mathcal{X}}(x^*)$ and $0 \in B^\top (Ax^*+B\bar{x}^*) +N_{\bar{\mathcal{X}}}(\bar{x}^*)$, which imply $(x^*, \bar{x}^*)$ is stationary for \eqref{eq: feasibility}.
\item Using the same case analysis as in part 2), we can see $\|z^k\|\rightarrow 0$ (along the subsequence converging to $z^*$). Thus taking limit on \eqref{eq: admm_stationary} completes the proof.
\end{enumerate}
\end{proof}
\subsection{Proof of Theorem \ref{thm:itercomplexity}}
\begin{proof}
   We use $T_k$ to denote an upper bound of the number of inner iterations of the $k$-th ADMM, which produces an $\epsilon$-stationary solution point of \eqref{eq: distributedOPF_relaxed} (see Definition \ref{def2}) . By Lemma 3 of \cite{sun2019two}, upon termination of ADMM, we can find $(x^k, \bar{x}^k, z^k, y^k)$ and corresponding $(d_1^k,d_2^k,d_3^k)$ such that
   \begin{subequations}\label{eq: admm stationary}
   \begin{align}
   		\|d_1^k\| = & \|\rho^kA^\top (B\bar{x}^{t-1}+z^{t-1}-B\bar{x}^t -z^t)\| \notag \\
   		       \leq & \rho^k \|A\|\|B\bar{x}^{t-1}+z^{t-1}-B\bar{x}^t -z^t\|,\\
   		\|d_2^k\| =&  \|\rho B^\top (z^t- z^{t-1})\|\leq \rho^k \|B\|\|z^{t-1} -z^t\|,\\	
   		\|d_3^k\| = & \frac{1}{2}\|z^{t-1} -z^t\|,
   \end{align}
    \end{subequations}
	where $1\leq t\leq T_k$ is some index during ADMM satisfying 
	\begin{align}\label{eq: bound ADMM res}
		& \|B\bar{x}^{t-1}-B\bar{x}^t\|+\|z^{t-1}-z^t\| \notag \\
		\leq & \sqrt{2}  (\|B\bar{x}^{t-1}-B\bar{x}^t\|^2+\|z^{t-1}-z^t\|^2)^{1/2} \notag \\
		\leq & \sqrt{2} \left(\frac{2(\overline{L}-\underline{L})}{\beta^k T_k} \right)^{1/2}
		= 2\left(\frac{\overline{L}-\underline{L}}{\beta^k T_k}\right)^{1/2}.
	\end{align}
	Recall $\rho^k = 2\beta^k = 2\beta^0c^k\geq 2\beta^0$. So \eqref{eq: admm stationary} and \eqref{eq: bound ADMM res} give
    \begin{equation}
     \max\{\|d_1^k\|,\|d_2^k\|, \|d_3^k\|\}\leq \rho^k \tau \left(\frac{\overline{L}-\underline{L}}{\beta^k T_k}\right)^{1/2},
    \end{equation}
    where $\tau = \max\{2\|A\|, 2\|B\|, 1/(2\beta^0)\}$.
    It is sufficient to find a $T_k$ with 
    \begin{equation}\label{eq: admm iter}
    	\frac{4\beta^0 (\overline{L}-\underline{L})\tau^2 c^k}{\epsilon^2} \leq T_k \leq \frac{4\beta^0 (\overline{L}-\underline{L})\tau^2 c^k}{\epsilon^2}+1,
    \end{equation}
    in order to get $\max\{\|d_1^k\|,\|d_2^k\|, \|d_3^k\|\}\leq \epsilon$. Let $K$ denote the number of outer-level ALM iterations. Then the total number of inner iterations is bounded by 
    \begin{equation}
    \sum_{k=1}^K T_k \leq \left\lceil \left(\frac{4\beta^0(\overline{L}-\underline{L})\tau^2 c}{c-1}\right)\left(\frac{c^{K}-1}{\epsilon^2}\right)\right\rceil +K.
    \end{equation} 
    
    Now it remains to bound the number of outer-level iterations $K$. Notice that the dual residuals $\|d_1^K\|$ and $\|d_2^K\|$ are controlled by $\epsilon$ at the end of the inner-level ADMM, so it is sufficient to ensure the primal residual $\|Ax^K+B\bar{x}^K\|\leq \epsilon$. By Theorem 3 of \cite{sun2019two}, 
    \begin{align*}
        \|Ax^K+Bx^K\|^2\leq \frac{2(\overline{L}-c(x^K)+ \langle \lambda^k, Ax^K+Bx^K\rangle)}{\beta^K},
    \end{align*}
    and the claimed $K_1$ is sufficient to ensure $\|Ax^{K_1}+Bx^{K_1}\|\leq \epsilon$. Next we show the claimed bound on $K_2$. If $\|\lambda^k + \beta^k z^k\|\leq \Lambda$ for all outer index $k$, then we have 
    \begin{align}
    \|z^K\| = \frac{\|\lambda^K+\beta^Kz^K-\lambda^K\|}{\beta^K}\leq \frac{\Lambda+M}{\beta^0 c^K}.
    \end{align}
    By \eqref{eq: admm stationary} and \eqref{eq: bound ADMM res} , we have
    \begin{align}
    	\|Ax^K+B\bar{x}^K+z^K\| = \|d_3^K\| \leq & \left(\frac{\overline{L}-\underline{L}}{\beta^KT_K}\right)^{1/2}\leq \frac{\epsilon}{2\beta^0c^K\tau}, 	\notag 
    \end{align}
	where the last inequality is due to the lower bound of $T_K$ in \eqref{eq: admm iter}. It is straightforward to verify that the claimed $K_2$ ensures
    {\begin{align*}
        \|Ax^{K_2}+B\bar{x}^{K_2}\| \leq &\|Ax^{K_2}+B\bar{x}^{K_2}+z^{K_2}\|+\|z^{K_2}\|\\
         \leq & \epsilon/2 + \epsilon/2 = \epsilon.
    \end{align*}}
    This completes the proof.
\end{proof}

\end{document}